\documentclass[10pt,a4paper]{article}
\usepackage[left=3cm,right=3cm]{geometry}

\usepackage{latexsym}
\usepackage{cite}
\usepackage{amstext}
\usepackage{amssymb}
\usepackage{amsmath}

\usepackage{amscd}
\usepackage{graphicx}
\usepackage[utf8]{inputenc}
\usepackage{bm}
\usepackage[english]{babel}
\usepackage{tikz}
\usepackage{subcaption}

\usepackage{tabularx}
\usepackage{verbatim}
\usepackage{enumerate}
\usepackage{xcolor}
\usepackage{hyperref}
\usepackage{float}
\usepackage{caption}

\usepackage{appendix}
\usepackage{array}
\usepackage{footnote}
\usepackage{booktabs}
\usepackage{hhline}

\usepackage{algpseudocode}
\usepackage{algorithm}
\usepackage{wrapfig}
\usepackage{cleveref}
\usepackage{amsopn}

\usepackage{stmaryrd}
\usepackage{amsmath}
\usepackage{tabularx}

\usepackage{lipsum}
\usepackage{amsfonts}
\usepackage{epstopdf}
\ifpdf
  \DeclareGraphicsExtensions{.eps,.pdf,.png,.jpg}
\else
  \DeclareGraphicsExtensions{.eps}
\fi

\newcommand{\R}{\mathbb{R}}
\DeclareMathOperator*{\argmin}{arg\,min}

\numberwithin{equation}{section}

\graphicspath{{../../images/}}

\newcommand{\prox}{\text{prox}}
\newcommand{\N}{\mathbb{N}}

\newtheorem{remark}{Remark}
\newtheorem{definition}{Definition}
\newtheorem{lemma}{Lemma}
\newtheorem{theorem}{Theorem}
\newtheorem{corollary}{Corollary}
\newtheorem{proposition}{Proposition}

\newcounter{auxFootnote}
\title{Parameter-Free FISTA by Adaptive Restart and Backtracking}

\author {J.-F. Aujol\footnote{Univ. Bordeaux, Bordeaux INP, CNRS, IMB, UMR 5251, F-33400 Talence, France}
\and
L. Calatroni\footnote{CNRS, Université Côte d’Azur Inria, Laboratoire I3S, France}
\and 
Ch. Dossal\footnote{IMT, Univ. Toulouse, INSA Toulouse, Toulouse, France}\setcounter{auxFootnote}{\value{footnote}}
\and
H. Labarri\`ere\footnotemark[\value{auxFootnote}] 
\and 
A. Rondepierre\footnotemark[\value{auxFootnote}] \footnote{LAAS, Univ. Toulouse, CNRS, Toulouse, France}\\
\footnotesize{Jean-Francois.Aujol@math.u-bordeaux.fr,}\\
\footnotesize{calatroni@i3s.unice.fr,}\\
\footnotesize{\{Charles.Dossal,Hippolyte.Labarriere,Aude.Rondepierre\}@insa-toulouse.fr}}
\begin{document}

\maketitle

\begin{abstract}
We consider a combined restarting and adaptive backtracking strategy for the popular Fast Iterative
Shrinking-Thresholding Algorithm \cite{Beck2009} frequently employed for accelerating the convergence speed of large-scale structured convex optimization problems. Several variants of FISTA enjoy a provable linear convergence rate for the function values $F(x_n)$ of the form $\mathcal{O}( e^{-K\sqrt{\mu/L}~n})$ under the prior knowledge of problem \emph{conditioning}, i.e. of the ratio between the (\L ojasiewicz) parameter $\mu$ determining the growth  of the objective function and the Lipschitz constant $L$ of its smooth component. These parameters are nonetheless hard to estimate in many practical cases. Recent works address the problem by estimating either parameter via suitable adaptive strategies. In our work both parameters can be estimated at the same time by means of an algorithmic restarting scheme where, at each restart, a non-monotone estimation of $L$ is performed. For this scheme, theoretical convergence results are proved, showing that a $\mathcal{O}( e^{-K\sqrt{\mu/L}n})$ convergence speed can still be achieved along with quantitative estimates of the conditioning. The resulting Free-FISTA algorithm is therefore parameter-free. Several numerical results are reported to confirm the practical interest of its use in many exemplar problems.
\end{abstract}



\section{Introduction}  \label{sec:intro}

The Fast Iterative Soft-Thresholding Algorithm (FISTA) has been popularized in the work of Beck and Teboulle \cite{Beck2009} as an extension of previous works by Nesterov \cite{Nesterov1983,Nesterov2004} where improved $O(1/n^2)$ convergence rate was shown upon suitable extrapolation of the algorithmic iterates. In \cite{Nesterov2004}, such rate is shown to be optimal for the class of convex functions, outperforming the one of the classical Forward-Backward algorithm \cite{CombettesWajs2005}. In its vanilla form, FISTA is indeed an efficient strategy for computing solutions of convex optimization problems of the form
\begin{equation}  \label{eq:prob_comp}
	\min_{x\in \R^N}~F(x) := f(x) + h(x) ,
\end{equation}
where $F:\R^N\to \R \cup \left\{+\infty\right\}$ belongs to $\mathcal{H}_{L}$, the class of composite functions with $f$  convex and differentiable with $L$-Lipschitz gradient and $h$ convex, proper and lower semicontinuous (l.s.c.) with simple (i.e. easily computable) proximal operator. We also assume: 
$X^*:=\argmin_x~F(x)\neq \emptyset$.

Due to its wide use in many areas of signal/image processing, many extensions of FISTA enjoying monotonicity \cite{beck2009fast}, general extrapolation rules \cite{Attouch2018}, inexact  proximal point evaluations \cite{VillaSalzo2013}, variable metrics \cite{Bonettini2018} and improved $o(1/n^2)$ convergence rate \cite{Attouch2016} were proposed along with a large number of FISTA-type algorithms addressing specific features (e.g., FASTA \cite{GoldsteinFASTA2014}, Faster-FISTA \cite{Liang2022} to  name a few).
The question on the convergence of iterates of FISTA was solved in \cite{ChambolleDossal2015} whose results were then further investigated in several other papers, see, e.g., \cite{Liang2017,Liang2022}. The algorithmic convergence of FISTA relies on an upper bound on the algorithmic step-size, which depends on the inverse of the Lipschitz constant $L$. Practically, the estimation of $L$ may be  pessimistic and/or costly, which may result in unnecessary small step-size values. To avoid this, several backtracking strategies have been proposed based either on monotone (Armijo-type) \cite{Beck2009} or adaptive updates \cite{Scheinberg2014}. 

Interestingly, when the function $F$ satisfies additional growth assumptions such as strong convexity or quadratic growth, first-order methods may provide improved convergence rates. Under such hypotheses, Heavy-Ball type methods provide the fastest convergence rates\footnote{We call Heavy-Ball methods the schemes that are derived from the Heavy-Ball with friction system which includes Polyak's Heavy-Ball method \cite{Polyak1964}, Nesterov's accelerated gradient method for strongly convex functions \cite{Nesterov2004}, iPiasco \cite{ochs2015ipiasco} or V-FISTA \cite[Section 10.7.7]{Beck2017}}. Such methods rely on a constant-in-time inertial coefficient which is chosen according to $\kappa=\frac{\mu}{L}$ where $\mu>0$ is the parameter appearing in the growth condition. In fact, $\kappa$ is the inverse of the condition number and knowing its value is crucial for these methods to reach rates of the form $\mathcal{O}\left(e^{-K\sqrt{\kappa}n}\right)$ for some real constant $K>0$. We refer the reader to \cite[Table~2]{aujol2020convergence} for further details and comparisons. Note that in such a setting the Forward-Backward method guarantees in fact a decay of the error in $\mathcal{O}\left(e^{-\kappa n}\right)$ which is much slower since $\kappa\ll1$ in general. Different approaches requiring the explicit prior knowledge of both strong convexity parameters $\mu_f$ and $\mu_h$ of the functions in \eqref{eq:prob_comp} have been studied in \cite{ChambollePock2016,Calatroni-Chambolle-2019,florea2020generalized} and endowed with possible adaptive backtracking strategies.

In \cite{aujol2021fista} it has been shown that unlike Heavy-Ball methods, FISTA does not significantly benefit from growth-type assumptions. The presence of an inertial coefficient growing with the iterations amplifies the effect of inertia, so the scheme can generate oscillations when the function $F$ is sharp. From  a theoretical viewpoint, the decay of the error cannot be better than polynomial although the finite-time behavior of FISTA is close to the one of Heavy-Ball methods.
Restarting FISTA for functions satisfying some growth condition is a natural way of controlling inertia, which allows to accelerate the overall convergence. The main idea consists in reinitializing to zero the inertial coefficient based on some restarting condition. Elementary computations show that by restarting every $k^*$ iterations for some $k^*$ depending on $\sqrt{\kappa}$, the worst-case convergence improves to $\mathcal{O}\left(e^{-K\sqrt{\kappa}~n}\right)$ for some $K>0$ \cite{zhang2015restricted,Ferocq2016,Necoara2019}. Nonetheless, such restarting rule requires the knowledge of $\kappa$ and provides slower worst-case guarantees than Heavy-Ball methods. On the other hand, adaptive restarting techniques allow the adaptation of the inertial parameters to $F$ without requiring any knowledge on its geometry (apart from $L$). In \cite{Candes2015}, the authors propose heuristic restart rules based on rules involving the values of $F$ or $\nabla F$ at each iterate. These schemes are efficient in practice as they do not require any estimate of $\kappa$, but they do not enjoy any rigorous convergence rate. Fercoq and Qu introduce in \cite{fercoq2019adaptive} a restarting scheme achieving a fast exponential decay of the error when only a (possibly rough) estimate  of $\mu$ is available. In \cite{Alamo2019,alamo2019gradient,alamo2022restart}, Alamo et al. propose strategies ensuring linear convergence rates only using information on $F$ or the composite gradient mapping at each iterate. Roulet and d'Aspremont propose  in \cite{roulet2020sharpness} a restarting scheme based on a grid-search strategy providing a fast decay as well. Note that by restarting FISTA an estimate of the growth parameter can be done as shown by Aujol et al. in \cite{FISTArestart}, where fast linear convergence is shown. 

Adaptive methods exploiting the geometry of $F$ without knowing its growth parameter $\mu$ are useful in practice since estimating $\mu$ is generally difficult. In the same spirit, numerical schemes for strongly convex functions where the growth parameter is unknown are provided by Nesterov in \cite[Section~5.3]{Nesterov2013} and by Gonzaga and Karas in \cite{gonzaga2013fine}. In the case of strongly convex objectives, Lin and Xiao introduced in \cite{Lin2015} an algorithm  achieving a fast exponential decay of the error by automatically estimating both $L$ and $\mu$ at the same time.

\medskip

In this paper we consider a parameter-free FISTA algorithm (called Free-FISTA) with provable accelerated linear convergence rates of the form $\mathcal{O}(e^{-K\sqrt{\kappa}n})$ for functions satisfying the quadratic growth condition:
\begin{equation}  \label{eq:growth}
    (\exists \mu>0)\quad\text{s.t.}\quad (\forall x\in\mathbb{R}^N)\quad\frac{\mu}{2}d(x,X^*)^2\leq F(x)-F^*,
\end{equation}
assuming that both the growth parameter $\mu>0$ and the Lipschitz smoothness parameter $L>0$ of $\nabla f$ are unknown. By a suitable combination of existing previous work combining an adaptive restarting strategy for the estimation of $\mu$ \cite{FISTArestart} and a non-monotone estimation of $L$ performed via adaptive backtracking at each restart \cite{Scheinberg2014,Calatroni-Chambolle-2019}, Free-FISTA adapts its parameters to the local geometry of the functional $F$, thus resulting in an effective performance on several exemplar problems in signal and image processing. The proposed strategy relies on an estimate ${\kappa}_j$ of $\kappa$ which is rigorously showed to provide a restarting rule that guarantees fast convergence. 

\section{Preliminaries and notations} \label{sec:preliminaries}

We are interested in solving the convex, non-smooth composite optimization problem \eqref{eq:prob_comp} under the following assumptions:
\begin{itemize}
    \item The function $f: \mathbb{R}^N\to\mathbb{R}_+$ is convex, differentiable with $L$-Lipschitz gradient:
    \begin{equation*}
    (\exists L\geq 0)\quad (\forall x,y\in \mathbb{R}^N)\quad \|\nabla f(x)-\nabla f(y)\|\leq L\|x-y\|.
    \end{equation*}
    \item The function $h:\mathbb{R}^N\to \mathbb{R}_+ \cup \left\{ +\infty\right\}$ is proper, l.s.c. and convex. Its proximal operator will be denoted by:
    \begin{equation}  \label{eq:prox_operator}
    \prox_{h}(z) = \argmin_{w\in \mathbb{R}^N}~h(w) + \frac{1}{2}\|w-z\|^2,\quad z\in\mathbb{R}^N.
    \end{equation}
\end{itemize}

For this class of functions a classical minimization algorithm is the Forward-Backward algorithm (FB) whose iterations are described by:
\begin{equation*}
  x_{k+1}=\prox_{\tau h}(x_k-\tau \nabla f(x_k)),\quad \tau\in \left(0,\frac{2}{L}\right).
\end{equation*}
To define in a compact way the Forward-Backward iteration performed on $y\in\mathbb{R}^N$ with a step-size $\tau>0$, we will use the notation
$T_\tau(y) = \prox_{\tau h}(y-\tau\nabla f(y)).$
while for assessing optimality via a suitable stopping criterion, we will consider a condition of the form $0\in\partial F(y)$, or, equivalently, $g_{\tau}(y)=0$  with the composite gradient mapping being defined by:
\begin{equation*}
g_{\tau}(y): = \frac{y-T_{\tau}(y)}{\tau}=\frac{1}{\tau}\left(y- \text{prox}_{\tau h}\left(y-\tau\nabla f(y)\right)\right),\quad y\in\mathbb{R}^N.
\end{equation*}
This last formulation is convenient for defining an approximate solution to the composite problem, and thus to deduce a tractable stopping criterion:
\begin{definition}[$\varepsilon$-solution]
Let $\varepsilon>0$ and $\tau>0$. An iterate $y\in\R^N$ is said to be an $\varepsilon$-solution of the problem
\eqref{eq:prob_comp} if: $\|g_\tau(y)\| \leqslant \varepsilon$.
\end{definition}
Given an estimation $\hat L>0$ of $L$ and a tolerance $\varepsilon>0$, the exit condition considered will then read $\|g_{1/\hat{L}}(y)\|\leqslant\varepsilon$.
As a shorthand notation, we also define the class of functions satisfying  \eqref{eq:growth}:
\begin{definition}[Functions with quadratic growth, $\mathcal{G}_\mu^2$]
Let $F:\R^N\rightarrow \R\cup \{+\infty\}$ be a proper l.s.c. convex function with $X^*:=\argmin~F \neq \emptyset$. Let $F^*:=\inf F$. The function $F$ satisfies a quadratic growth condition $\mathcal{G}_\mu^2$ for some $\mu>0$ if:
\begin{equation}
(\forall x\in\R^N),\qquad \frac{\mu}{2}d(x,X^*)^2\leqslant F(x)-F^*.
\label{eq:Lojasiewicz_dis}
\end{equation}
\end{definition}

Condition \eqref{eq:Lojasiewicz_dis} can be seen as a relaxation of strong convexity. As shown in \cite{Bolte2017,garrigos2017convergence} in a convex setting such condition is equivalent to a global \L ojasiewicz property with an exponent $\frac{1}{2}$. In particular, the following lemma states an implication that is required in the later sections.

\begin{lemma}
\label{lem:Loja}
Let $F:\R^N\rightarrow \R\cup \{+\infty\}$ be a proper, l.s.c. and convex function with a non-empty set of minimizers $X^*$. Let $F^*=\inf F$. If $F$ satisfies $\mathcal{G}^2_\mu$ for some $\mu>0$, then $F$ has a global \L{}ojasiewicz property with an exponent $\frac{1}{2}$:
\begin{equation*}
(\forall x\in \R^N),\quad \frac{\mu}{2}\left(F(x)-F^*\right)\leqslant d(0,\partial F(x))^2.
\end{equation*}
\end{lemma}


\section{Free-FISTA}  \label{sec:free_fista}
In this paper we propose a parameter-free restart algorithm based on the original FISTA scheme proposed by Beck and Teboulle in \cite{beck2009fast}:
\begin{equation*}
 y_k=x_k+ \frac{t_k-1}{t_{k+1}}(x_k - x_{k-1}),\quad
x_{k+1}=\prox_{\tau h}(y_k-\tau \nabla f(y_k)),
\end{equation*}
where the sequence $(t_k)_{k\in\N}$ is recursively defined by: $t_1=1$ and $t_{k+1} =(1+\sqrt{1+4t_k^2})/{2}$. For the class of convex composite functions, the convergence rate of the method is given by \cite{Nesterov1983,beck2009fast}:
\begin{equation*}
(\forall k\in \N),\quad F(x_k)-F^* \leqslant \frac{2L\|x_0-x^*\|^2}{(k+1)^2}.
\end{equation*}
When $L$ is available, a classical strategy introduced in \cite{nesterov27method} is to restart the algorithm at regular intervals. Necoara and al. \cite{Necoara2019} propose an optimized restart scheme, proving that restarting Nesterov accelerated gradient every $\lfloor 2e\sqrt{\frac{L}{\mu}}\rfloor$ iterations ensures that $F(x_k)-F^*=\mathcal{O}\left(e^{-\frac{1}{e}\sqrt{\frac{\mu}{L}}k}\right)$ for the class of $\mu$-strongly convex functions. This restart scheme and its convergence analysis can be extended to composite functions satisfying some quadratic growth condition $\mathcal{G}^2_\mu$ \cite{Candes2015,Necoara2019}.

In this paper we consider the case when both the Lipschitz constant $L$ and the growth parameter $\mu$ are unknown. The first main ingredient of our parameter-free FISTA algorithm is the use of an adaptive backtracking strategy used at each restart to provide a non-monotone estimation of the local Lipschitz constant $L$. More precisely, we propose a backtracking variant of FISTA (FISTA-BT), widely inspired by the one proposed in \cite{Calatroni-Chambolle-2019} and described in Section \ref{sct:backtracking}. The second main ingredient is an adaptative restarting approach, described in Section \ref{sct:adaptive:restart}, taking advantage of the local estimation of the geometry of $F$ (via online estimations of the parameter $\kappa=\frac{\mu}{L}$) for avoiding oscillations due to inertia. The main steps of Free-FISTA are the following: at each restart, given a current iterate $r_{j-1}$, a fixed number of iterations $n_{j-1}$ and a current estimation $L_{j-1}^+$ of the Lipschitz constant $L$,
\begin{enumerate}
\item Compute $r_j$ a new iterate and $L_j$ a new estimation of $L$ by performing $n_{j-1}$ iterations of FISTA-BT algorithm parameterized by the estimate $L_{j-1}^+$.
\item Compute an estimation $\kappa_j$ of the geometric parameter $\kappa=\frac{\mu}{L}$.
\item Update the number $n_j$ of iterations of FISTA-BT for the next restart loop. It depends on $n_{j-1}$ and on $\kappa_j$.
\end{enumerate}
The whole algorithm is carefully described in Section \ref{sct:FreeFISTA} and its convergence is proven.
All technical proofs are reported in a dedicated \Cref{sec:proof}.

\subsection{Adaptive backtracking}   \label{sct:backtracking}

In order to provide at each restart of  Free-FISTA an estimation of $L$ adapted to the current estimate of the growth parameter, we describe in the following an instance of FISTA endowed with non-monotone backtracking previously considered, e.g., in \cite[Algorithm 2]{Scheinberg2014} and \cite[Algorithm 2]{Calatroni-Chambolle-2019} with $\mu=0$.
Differently from standard approaches following an Armijo-type (i.e. monotone) backtracking rule \cite{Beck2009}, the use of a non-monotone strategy further allows for a local decreasing of the estimated valued $\hat{L}$ of $L$ (equivalently, an increasing of $\tau$ w.r.t. to the optimal $1/L$)  in the neighborhoods of ``flat" points of the function $f$ (i.e. where $L$ is small), thus improving practical performances.

Following \cite{Calatroni-Chambolle-2019}, the proposed adaptive backtracking strategy is derived from the classical descent condition holding for FISTA at ${x}^+:= T_{\tau}(x)$ with $x\in\R^N$, which reads: for any $y\in \R^N$,
\begin{equation}   \label{condition:backtr}
F(x^+) + \frac{\|y-x^+\|^2}{2\tau} + \left(  \frac{\|x^+-x\|^2}{2\tau} - D_f(x^+,x) \right) \leq F(y) + \frac{\|y-x^+\|^2}{2\tau},
\end{equation}
which is defined in terms of the Bregman divergence $D_f:\mathbb{R}^N\times \mathbb{R}^N\to\mathbb{R}_+$ associated to $f$ and defined by: 
$
D_f(x,y) = f(x) - f(y) - \langle \nabla f(y),x-y\rangle.
$
Choosing $y=x$ in \eqref{condition:backtr}, the descent of $F$ between two iterates $x$ and $x^+=T_\tau(x)$ is at least of:
\begin{equation}  \label{condition:bregman}
F(x^+)-F(x) \leqslant  -\frac{\|x^+-x\|^2}{2\tau}, \quad \text{provided that}\quad 
D_f(x^+,x)\leq \frac{\|{x}^+-{x}\|^2}{2\tau}.   
\end{equation}
This last condition is true whenever $0<\tau\leq 1/L$. When only a local estimate $L_k$ of $L$ is available, the idea is to enforce \eqref{condition:bregman} by applying a backtracking strategy t by $\tau_k=\frac{1}{L_k}$: testing a tentative step-size $\tau_{k}=\tau_{k-1}/\delta$ with $\delta\in (0,1)$ greater than the one $\tau_{k-1}$ considered at the previous iteration, decrease the step $\tau_{k}$ by a factor $\rho\in (0,1)$ as long as condition \eqref{condition:bregman} is not satisfied. This condition can be rewritten as
$\frac{2D_f(x^+,x)}{\|x^+-x\|^2}>\frac{\rho}{\tau_k} = \rho L_k$, 
where $\tau_k/\rho$ denotes the last step before acceptance. Note that by the condition above, for all $k\geq 0$ there holds:
\begin{equation}    \label{eq:Lk_rho}
\tau_k\geq \frac{\rho}{L} \qquad \Leftrightarrow \qquad L_k \leq \frac{L}{\rho},
\end{equation}
which can be used to get the desired convergence result.

The algorithm \texttt{FISTA\_adaBT} is reported in  \Cref{alg:FISTABT}.
\begin{algorithm}[h]
  \caption{FISTA + adaptive backtracking, \texttt{FISTA\_adaBT}($x^0,n,L_0, L_{min};\rho, \delta$)}
\label{alg:FISTABT}
  \begin{algorithmic}
    
\Statex {\textbf{Initializations}: $\tau_0=1/L_0$, $\rho\in(0,1), \delta\in(0,1]$, $x_{-1}=x_0 \in \mathcal{X}$, $t_0=1$, $L_{min}$ sufficiently small.}

    \vspace{0.1cm}

    \For  {$k=0,1,\ldots,n$}   
 \begin{equation}  \label{alg:increase}
 \tau^0_{k+1}=\min\left\{\frac{\tau_{k}}{\delta},\frac{1}{L_{min}}\right\};
 \end{equation}  
\qquad$i=0;$\\
\qquad{\textbf{repeat}} 
{
 {
 \begin{align}
 \tau_{k+1} & = \rho^{i}~\tau^0_{k+1}; \notag \\
 t_{k+1}& = \frac{1 + \sqrt{1 + 4\frac{\tau_{k}}{\tau_{k+1}}t^2_k}} {2}; \label{update:tk}\\   
\beta_{k+1} & = \frac{t_k-1}{t_{k+1}}; \notag \\
 y_{k+1} & = x_{k} + \beta_{k+1}(x_{k}-x_{k-1});  \notag \\
x_{k+1} & = \text{prox}_{\tau_{k+1} h} (y_{k+1} - \tau_{k+1}\nabla f(y_{k+1})); \notag \\
i&= i+1 ;\notag
\end{align}

}

}
 \textbf{until ~$D_f(x_{k+1},y_{k+1})\leq \|x_{k+1}-y_{k+1}\|^2/2\tau_{k+1}$}


    \EndFor

\noindent \textbf{Return} $(x_{k+1},L_{k+1}=\frac{1}{\tau_{k+1}})$
  \end{algorithmic}
\end{algorithm}
The parameter $L_{min}>0$ provides a lower bound of the estimated Lipschitz constants at any $k$, i.e $L_k=\frac{1}{\tau_k}\geqslant L_{min}$. This property will be needed to prove the theoretical asymptotic convergence rate of the global restarting scheme. Such parameter has to satisfy the condition $L_{min}<L$. However, since this value should be taken as small as possible this condition is not restrictive and it practically does not affect the choice \eqref{alg:increase}. We observe that whenever $\delta<1$, the increasing of the algorithmic step-size is attempted at each outer iteration of \Cref{alg:FISTABT}, while, when $\delta=1$, the same value $\tau_k$ estimated at the previous iterations is used. In both cases, a standard Armijo-type backtracking is then run to adjust possible over-estimations.



Convergence of \Cref{alg:FISTABT} is stated in the following Theorem, which is a special case of \cite[Theorem 4.6]{Calatroni-Chambolle-2019} suited for the particular case $\mu=0$ (no strong-convexity). 

\begin{theorem}[Convergence of \Cref{alg:FISTABT} \cite{Calatroni-Chambolle-2019}]  \label{theo:convergence}
Let $n\in \N$. The sequence $(x_k)_{k=0,\dots, n}$ generated by the  Algorithm \ref{alg:FISTABT} satisfies for all $k=0,\ldots,n$:
\begin{equation}  \label{convergence:functional2}
F(x_{k+1})-F^*\leq \frac{2 \bar{L}_{k+1}}{(k+1)^2}  \|x_0-x^*\|^2,
\end{equation}
where, by setting $L_i:=1/\tau_i$ the quantity $\bar{L}_{k+1}$ is defined by:
\begin{equation}   \label{eq:local_estimate}
	\bar{L}_{k+1} :=   \left( \frac{1}{\frac{1}{k+1} \sum_{i=1}^{k+1}  \frac{1}{\sqrt{L_i}}   } \right)^2.
\end{equation}

\end{theorem}

The (harmonic) average appearing in  \eqref{convergence:functional2} depends only on the estimates of $L$ performed along the iterations of \Cref{alg:FISTABT}. In particular, it does not depend on the unknown value of the Lipschitz constant $L$. However, recalling \eqref{eq:Lk_rho}, we have for all $k=1,\dots,n$, $\rho\bar{L}_{k+1}\leqslant L$, hence the following bound:
\begin{equation}  \label{eq:conv_FISTA_BT}
\frac{2 \bar{L}_{k+1}}{(k+1)^2}\leq \frac{2L}{\rho(k+1)^2}
\end{equation}
which, plugged in \eqref{convergence:functional2},  entails the well-known convergence rate for  FISTA endowed with Armijo-type backtracking showed, e.g., in \cite{Beck2009}.

\begin{remark}
Regarding the choice of the extrapolation rule \eqref{update:tk}, we remark that in  \cite{FISTArestart} a different update based on \cite{ChambolleDossal2015} was considered to guarantee the convergence of the iterates of the resulting FISTA scheme. Since the convergence result in Theorem \ref{theo:convergence} cannot be adapted to this different choice in a straightforward manner, we consider in this work a Nesterov-type update, inspired by previous work \cite{Scheinberg2014,Calatroni-Chambolle-2019}.
\end{remark}

 We can now state the main proposition (whose proof is detailed in \Cref{sec:prop1}) which will be used in the following to formulate the proposed adaptive restarting strategy described in \Cref{sct:adaptive:restart}:
\begin{proposition}
\label{prop:Nesterov_up}
Let $F$ be a function satisfying $\mathcal{H}_L$ and $\mathcal{G}^2_\mu$ for some $L>0$ and $\mu>0$. If $L_{min}\in[0,L)$, then for any fixed $n\in \mathbb{N}^*$, the sequence $(x_k)_{k=0\dots n}$ provided by \Cref{alg:FISTABT} satisfies for all $k\in \N$:
\begin{align}
\text{(i)}&~\ F(x_{k+1})-F^*\leqslant\frac{4L}{\rho\mu (k+1)^2}\left(F(x_0)-F^*\right),\label{eq:Nesterov_up1}\\
\text{(ii)}&~ F(x_{k+1})\leqslant F(x_0),\label{eq:Nesterov_up2}
\end{align}
\end{proposition}


\subsection{Adaptive restarting}\label{sct:adaptive:restart}

Having provided an estimate of $L$ after one algorithmic restart $j\geq 1$, intuitively, let us now describe the strategy of Free-FISTA. The structure of the algorithm relies on two main ingredients: a tractable stopping criterion suitable to cope with the hypothesis that the Lipschitz constant $L$ is not available, and a strategy to approximate the unknown value of the conditioning parameter $\kappa=\frac{\mu}{L}$ by a sequence $(\kappa_j)_j$ whose values will be needed to define the number $n_j$ of inner FISTA-BT iterations to be performed at each restart.

\subsubsection{A tractable stopping criterion}

Let $\varepsilon>0$ be the expected accuracy and $(r_j,L_j)$ be the $j-th$ output of \Cref{alg:FISTABT} for $n_{j-1}$ iterations at the $j-th$ restart.
When the Lipschitz constant $L$ is available, the notion of $\varepsilon$-solution can be seen as a good stopping criterion for an algorithm solving the composite optimization problem for three reasons: first it is numerically quantifiable. Secondly controlling the norm of the composite gradient mapping is roughly equivalent to having a control on the values of the objective function. Lastly, it will enable to analyze and compare algorithms in terms of the number of iterations needed to reach the accuracy $\varepsilon$.

\begin{algorithm}[H]
	\caption{Forward-Backward step with Armijo-backtracking, \texttt{\texttt{FB\_BT}}($r,L_0;\rho$)}
	\label{alg:FB_BT_up}
	\begin{algorithmic}
		\Statex {\textbf{Require}: $r\in\R^N$, $L_0>0$, $\rho\in(0,1)$.}

  \vspace{0.1cm}

$i=0$
  
\Repeat 

$\tau=\frac{\rho^i}{L_0}$

${r}^+ = T_\tau (r)$

$i=i+1$

\Until{$D_f(r^+,r)\leqslant \|r^+-r\|^2/2\tau$}\\
\textbf{Return} $r^+,~L^+=\frac{L_0}{\rho^{i-1}}$
\end{algorithmic}
\end{algorithm}

When only estimations $L_j$ of $L$ are available at each restart, there is no guarantee that the condition  $\|g_{1/L_j}(r_j)\|\leqslant \varepsilon$ will enable to control the values of the objective functions. To get a tractable stopping criterion, we propose to add a Forward-Backward step with Armijo backtracking before the next restart. Such an algorithm, denoted by \texttt{FB\_BT}, is detailed in \Cref{alg:FB_BT_up}. This extra step ensures that the following condition holds for all $j\geqslant1$:
\begin{equation}  \label{eq:breg_FB_BT}
D_f({r}_j^+,r_j)\leqslant \frac{{L}_j^+}{2}\|{r}_j^+-r_j\|^2,
\end{equation}
where $({r}_j^+=T_{1/{L}_j^+}(r_j),{L}_j^+)$ denote the outputs of \Cref{alg:FB_BT_up}, and $
g_{1/L_j^+}(r_j) = L_j^+(r_j-r_j^+)$
with, by construction: $L_j^+\geqslant L_j$. Note that the computational cost of the composite gradient mapping $g_{1/L_j^+}(r_j)$ is therefore very low. The stopping criterion of Free-FISTA thus reads:
\begin{equation}
\|g_{1/L_j^+}(r_j)\|\leqslant \varepsilon.\label{exit:condition}
\end{equation}
The condition \eqref{exit:condition} is a ``good" stopping criterion in the sense that it enables to control the values of the objective function along the iterations. Our analysis relies on the following Lemma whose proof is detailed in \Cref{sec:lem4}:
\begin{lemma}
\label{lem:cor}
Let $F$ be a function satisfying $\mathcal{H}_L$ and $\mathcal{G}^2_\mu$ for some $L>0$ and $\mu>0$. Then for all $x\in\R^N$ and $\tau>0$ we have:
\begin{equation*}
F(T_\tau (x))-F^*\leqslant \frac{2(1+L\tau)^2}{\mu}\|g_\tau(x)\|^2.
\end{equation*}
\end{lemma}
Applying Lemma \ref{lem:cor} to the iterate $r_j$, we get:
\begin{equation*}
F({r}_j^+)-F^*\leqslant\frac{2(1+L/{L}_j^+)^2\varepsilon^2}{\mu},
\end{equation*}
where, importantly, does not require the computation of $F^*$. 
In addition, remembering that the parameter $L_{min}\in\left(0,L\right)$ from \Cref{alg:FISTABT} provides a lower bound on the estimates $L_j$ and that $L_j^+\geqslant L_j$, we necessarily have: $L_j^+ \geqslant L_{min}$ and thus:
\begin{equation*}
F({r}_j^+)-F^*\leqslant\frac{2(1+L/L_{min})^2\varepsilon^2}{\mu}.
\end{equation*}

\begin{remark}
An alternative choice for ${L}_j$ following from \eqref{convergence:functional2} is ${L}_j=\bar{L}_j$ with
\[
\bar{L}_j = \left( \frac{1}{\frac{1}{n_{j-1}} \sum_{k=1}^{n_{j-1}}  \frac{1}{\sqrt{L_k}}   } \right)^2
\]
being the average \eqref{eq:local_estimate} estimated at the $j$-the restart. Nonetheless, we prefer ${L}_j=\frac{1}{\tau_{n_{j-1}}}\leqslant \frac{L}{\rho}$, as  the last estimation of $L$ at the $j$-th restart approximates the local smoothness of the functional. Moreover, its value is in general smaller than the value $\bar{L}_j$, which, when used for the next call of \Cref{alg:FISTABT} is expected to require fewer adjustments, thus improving the overall efficiency.
\end{remark}

%

\subsubsection{Estimating the geometric paramater $\kappa$}
Once the stopping criterion is well defined, the next issue is to determine the number of FISTA-BT iterations to perform at each restart. The global principle of our restart scheme is as follows: at the $j$-th restart,
\begin{itemize}
\item Compute $(r_j,L_j) = \mbox{FISTA\_adaBT}(r_{j-1}^+, n_{j-1},L_{j-1}^+,L_{min};\rho,\delta)$ where $r_j$ is the iterate computed after $n_{j-1}$ iterations of FISTA\_adaBT and $L_j$ the associated estimate of the Lipschitz constant $L$.
\item Perform an extra step of backtracking Forward-Backward: 
$$(r_j^+,L_j^+)=\mbox{FB\_BT}(r_j,L_j;\rho).$$
\item Update the number of iterations $n_j$ for the next restart.
\end{itemize}
Inspired by \cite{FISTArestart}, the update of the number $n_j$ of iterations relies on the estimation of the inverse $\kappa=\frac{\mu}{L}$ of the conditioning  at each restart loop by comparing the values $F(r_j)-F^*$ and $F(r_{j-1})-F^*$ at each restart $j$. More precisely, applying the first claim of \Cref{prop:Nesterov_up} at the $j$-th restart, we have: for all $j\in\N^*$\small
\begin{eqnarray*}
F(r_j)-F^*&\leqslant&\frac{4L}{\rho\mu (n_{j-1}+1)^2}\left(F( r_{j-1}^+)-F^*\right)
\leqslant\frac{4L}{\rho\mu (n_{j-1}+1)^2}\left(F(r_{j-1})-F^*\right),
\label{eq:Nesterov_rate}
\end{eqnarray*}\normalsize
observing that by the property \eqref{eq:breg_FB_BT}, we have: $F(r_j^+)\leqslant F(r_j)$ as explained in \Cref{sct:backtracking}. 
We thus deduce:
\begin{equation}
(\forall j\in\N^*),\quad \kappa\leqslant\frac{4}{\rho(n_{j-1}+1)^2}\frac{F(r_{j-1})-F^*}{F(r_j)-F^*}.\label{eq:ineq_muL}
\end{equation}
Since $F^*$ is often not known in practice and noticing that the application $u\mapsto \frac{F(r_{j-1})-u}{F(r_j)-u}$ is non decreasing on $[F^*,F(r_j)]$ (since $F(r_j)\leqslant F(r_{j-1})$), we deduce: 
\begin{equation*}
(\forall j\in\N^*),\quad \kappa\leqslant\frac{4}{\rho(n_{j-1}+1)^2}\frac{F(r_{j-1})-F(r_{j+1})}{F(r_j)-F(r_{j+1})}.
\end{equation*}
Using such inequality, it is thus possible to get a sequence $\left( \kappa_j\right)_j$ estimating $\kappa$ at each restart  $j\geqslant2$ by comparing $F(r_{j-1})-F(r_{j})$ and $F(r_{j-2})-F(r_{j})$ by defining:
\begin{equation}\label{eq:mu_j}
(\forall j\geqslant 2),\quad  \kappa_j:=\min\limits_{\underset{i<j}{i\in\N^*}}~\frac{4}{\rho(n_{i-1}+1)^2}\frac{F(r_{i-1})-F(r_{j})}{F(r_i)-F(r_{j})}.\end{equation}
By construction the sequence $(\kappa_j)_{j\in \N}$ is non-increasing along the iterations :
\begin{lemma}
Let $F$ be a function satisfying $\mathcal{H}_L$ and $\mathcal{G}^2_\mu$ for some $L>0$ and $\mu>0$. Then the sequence $( \kappa_j)_{j\geqslant2}$ defined by \eqref{eq:mu_j} satisfies
\begin{equation}(\forall j\geqslant 2),\quad \kappa_j\geqslant \kappa_{j+1}>\kappa.\end{equation}
\label{lem:mu_j}
\end{lemma}

\subsection{Free-FISTA: structure and convergence results} \label{sct:FreeFISTA}
Free-FISTA is detailed in \Cref{alg:FISTArestartBT_up}. Note that the `free' dependence on parameters stressed here relates to the two smoothness and growth parameters, $L$ and $\mu$, respectively. The  hyperaparameters $\rho\in(0,1)$, $\delta\in(0,1]$, $L_{min}>0$ required by Free-FISTA to perform adaptive backtracking  and to assess the expected precision ($0<\varepsilon\ll 1$) do not affect its convergence properties.

\begin{algorithm}[h!]
	\caption{Free-FISTA: parameter-free FISTA with adaptive backtracking and restart}
	\label{alg:FISTArestartBT_up}
	\begin{algorithmic}
		\Statex {\textbf{require}: $r_0\in\R^N$, $j=1$, $L_0>0$, $L_{min}>0$, $\rho\in(0,1)$, $\delta\in(0,1]$, $0<\varepsilon\ll 1$}

$n_0 = \lfloor 2C \rfloor$

$(r_1, {L}_1) =$\texttt{FISTA\_adaBT}$(r_0,n_0,L_0,L_{min};\rho,\delta)$

$n_1 = \lfloor 2C \rfloor$

$({r}_1^+,{L}_1^+)=\text{\texttt{FB\_BT}}(r_1,{L}_1;\rho)$

\Repeat

\vspace{0.1cm}

	$j=j+1$
	 
	 $(r_j, {L}_j) =$\texttt{FISTA\_adaBT}$({r}_{j-1}^+,n_{j-1},{L}_{j-1}^+,L_{min};\rho,\delta)$
	 
	${\kappa}_j = \min_{i<j } ~ \frac{4}{\rho(n_{i-1}+1)^2} \frac{F(r_{i-1}) - F(r_j)}{F(r_i)- F(r_j)} $
	
	\If {$n_{j-1} \leq C \sqrt{\frac{1}{{\kappa}_j}}$}
	
	$n_j=2n_{j-1}$
	
	\Else
	
	$n_j=n_{j-1}$
	
	\EndIf
	 	
	$({r}_j^+,L_j^+)=\text{\texttt{FB\_BT}}(r_j,{L}_j;\rho)$
	 
\Until{$\| g_{1/L_j^+}(r_j) \|  \leqslant \varepsilon $}\\
\Return $r={r}_j^+$

\end{algorithmic}
\end{algorithm}

To summarize, Free-FISTA \Cref{alg:FISTArestartBT_up} relies on a few sequences:
\begin{itemize}
	\item the sequence $(r_j)_{j\in\N}$ corresponds to the (outer/global) iterates. For all $j>0$, $r_j$ is the output of the $j$-th execution of \Cref{alg:FISTABT} after one extra application of \Cref{alg:FB_BT_up}.
	\item the sequence $(n_j)_{j\in\N}$ refers to the number of estimated iterations of \Cref{alg:FISTABT} to be performed at the $j$-th restart. For all $j\geqslant0$ we thus have:
	$$
	(r_{j+1},{L}_j)=\text{\texttt{FISTA\_adaBT}}({r}_j^+,n_j;{L}_j^+,L_{min};\rho,\delta),
	$$
 where $(r_j^+,L_j^+)$ is obtained after an extra Forward-Backward step with backtracking applied to $(r_j,L_j)$.
	\item the sequence $({L}_j)_j$ estimating $L$ at each restart. 
	\item the sequence $(\kappa_j)_{j\geqslant2}$ estimating  at each restart the true problem conditioning $\kappa=\mu/L$ by comparing the cost function $F$ at three different iteration points.
\end{itemize}

Let us finally explain our strategy to update the number $n_j$ of iterations required by \Cref{alg:FISTABT} at the $j$-th restart. Once an estimate $ \kappa_j$ is computed, the strategy performed by Free-FISTA consists in updating $n_j$ using a doubling condition that depends on a parameter $C>0$ to be defined:
\begin{equation}
n_{j-1} \leq C \sqrt{\frac{1}{{\kappa}_j}}\label{doubling:cd}
\end{equation}
Thus, Free-FISTA checks whether such condition is fulfilled: if it holds true, then $n_{j-1}$ is considered too small and doubled so that $n_j=2n_{j-1}$. Otherwise, the number of iterations is kept unchanged. By construction, the sequence $(n_j)_{j\in\N}$ is non-decreasing, and satisfies the following lemma.
\begin{lemma}
\label{lem:n_j}
Let $F$ be a function satisfying $\mathcal{H}_L$ and $\mathcal{G}^2_\mu$ for some $L>0$ and $\mu>0$. Then the sequence $(n_j)_{j\in\N}$ provided by \Cref{alg:FISTArestartBT_up} satisfies 
\begin{equation*}
\forall j\in\N,\quad n_j\leqslant2C\sqrt{\frac{1}{\kappa}}.
\end{equation*}
\end{lemma}

Note that for all $j\geqslant 2$, the number of iterations $n_j$ is defined according to $n_{j-1}$, $\kappa_j$ and the predefined parameter $C>0$. The proof of Lemma \ref{lem:n_j} is straightforward by induction: first observe that $n_0 = \lfloor 2C\rfloor \leqslant 2C \leqslant2C\sqrt{\frac{1}{\kappa}}$. Assume that $n_{j-1}\leqslant 2C\sqrt{1/\kappa}$. By construction, either \eqref{doubling:cd} is satisfied and
$n_j=2n_{j-1} \leqslant 2C\sqrt{\frac{1}{\kappa_j}}\leqslant 2C\sqrt{\frac{1}{\kappa}},$
by monotonicity of $(\kappa_j)_{j\in \N}$ (see \Cref{lem:mu_j}), or \eqref{doubling:cd} is not satisfied, and $n_j=n_{j-1}\leqslant 2C\sqrt{1/\kappa}$ by assumption. 




We can now state the main convergence results of Free-FISTA.  Their proof can be found in \Cref{sec:proof_thm1} and \Cref{sec:proofcor1}, respectively.

\begin{theorem}
\label{thm:1}

Let $F$ be a function satisfying $\mathcal{H}_L$ and $\mathcal{G}^2_\mu$ for some $L>0$ and $\mu>0$. Let $(r_j)_{j\in\N}$ and $(n_j)_{j\in\N}$ be the sequences provided by \Cref{alg:FISTArestartBT_up} with parameters $C>4/\sqrt{\rho}$ and $\varepsilon>0$. Then, the number of iterations $1+\sum_{i=0}^jn_i$ required to guarantee $\|g_{1/{L}_j^+}(r_j)\|\leqslant\varepsilon$ is bounded and satisfies\footnotesize
\begin{equation*}
\sum_{i=0}^j n_i \leqslant \frac{4C}{\log\left(\frac{C^2\rho}{4}-1\right)}\sqrt{\frac{L}{\mu}}\left(2\log\left(\frac{C^2\rho}{4}-1\right)+\log\left(1+\frac{16}{C^2\rho-16}\frac{2L(F(r_0)-F^*)}{\rho\varepsilon^2}\right)\right).
\end{equation*}\normalsize
\end{theorem}
\begin{corollary}
Let $F$ be as above. If $C>4/\sqrt{\rho}$, $\varepsilon>0$ and $L_{min}\in(0,L)$, then the sequences $(r_j)_{j\in\N}$ and $(n_j)_{j\in\N}$ provided by \Cref{alg:FISTArestartBT_up} satisfy
\begin{equation*}
F({r}_j^+)-F^*= \mathcal{O}\left(e^{-\frac{\log\left(\frac{C^2\rho}{4}-1\right)}{4C}\sqrt{\kappa}\displaystyle\sum_{i=0}^jn_i}\right).
\end{equation*}
Moreover, the trajectory of total number of FISTA iterates has a finite length and the method converges to a minimizer $x^*\in X^*$.
\\
Specifically, if $C$ maximizes $\frac{\log\left(\frac{C^2\rho}{4}-1\right)}{4C}$, namely $C\approx6.38/\sqrt{\rho}$, then there exists $K>\frac{1}{12}$ such that the sequences $(r_j)_{j\in\N}$ and $(n_j)_{j\in\N}$ satisfy
\begin{equation*}
F({r}_j^+)-F^*= \mathcal{O}\left(e^{-\sqrt{\rho}K\sqrt{\kappa}\displaystyle\sum_{i=0}^jn_i}\right).
\end{equation*}
\label{cor:1}
\end{corollary}



\Cref{cor:1} states that the Free-FISTA algorithm \ref{alg:FISTArestartBT_up} provides asymptotically a fast exponential decay. This convergence rate is consistent with the one expected for functions $F$ satisfying $\mathcal{H}_L$ and $\mathcal{G}^2_\mu$ where both the parameters $L$ and $\mu$ are unknown a priori. Note that in this setting Forward-Backward algorithm provides a low exponential decay
The variation of Heavy-Ball method introduced in \cite{aujol2020convergence}, the FISTA restart scheme introduced in \cite{fercoq2019adaptive} and fixed restart of FISTA require to estimate the growth parameter to ensure a fast exponential decay.
FISTA algorithm has the same fast decay as Free-FISTA in finite time (see \cite{aujol2021fista}), but with a smaller constant.

\section{Numerical experiments}

In this section, we report several applications of the Free-FISTA Algorithm \ref{alg:FISTArestartBT_up} showing how an automatic estimation of the smoothness parameter $L$ and the growth parameter $\mu$ can be beneficial. The combined approach is compared with vanilla FISTA \cite{Beck2009}, FISTA with restart \cite{FISTArestart} and FISTA with adaptive backtracking (\Cref{alg:FISTABT}) \cite{Calatroni-Chambolle-2019}. The first two examples show the advantages of Free-FISTA in comparison with other schemes, while the last example highlights some existing limitations of restarting methods. The codes that generate the figures are available in the following GitHub repository:
\url{https://github.com/HippolyteLBRRR/Benchmarking_Free_FISTA.git}

\subsection{Logistic regression with $\ell^2$-$\ell^1$-regularization}
As a first example, we focus on a classification problem defined in terms of a given dictionary $A\in \mathbb{R}^{m\times n}$ and labels $b\in\left\{-1,1\right\}^m$. We consider the minimization problem:
\begin{equation}
\min_{x\in\mathbb{R}^n}~ F(x):=\underbrace{\frac{\lambda_1}{2\|A^Tb\|_\infty}\sum_{j=1}^m \log \left(1+e^{-b_ja_j^T x}\right)+\frac{\lambda_2}{2}\|x\|^2}_{:=f(x)}+\underbrace{\|x\|_1}_{:=h(x)},\label{eq:reglog}\end{equation}
where $a_j=(A_{i,j})_j\in\mathbb{R}^m$ is the $j$-th row of $A$, $\lambda_1>0$ and $\lambda_2>0$. By definition, the value $x^*$ minimizing $F$ is expected to satisfy $\mathbb{P}(b_i=1|a_i)=\frac{1}{1+e^{-a_i^T x^*}}$ for any $i\in\llbracket 1,m\rrbracket$.  Note that the $\ell^2$ term aims to smooth the objective function while the $\ell^1$ regularization  sparsfies the solution which helps preventing from overfitting. An upper estimation of $L$ can easily be computed:
\begin{equation} \label{eq:L_estimate_prob1}
\hat L=\frac{\lambda_1\|A^Tb\|^2}{8\|A^Tb\|_\infty}+\lambda_2,
\end{equation}
which may be large whenever $\|A^T b\|\gg 1$.
We note that the function $F$ satisfies the assumption $\mathcal{G}^2_\mu$ for some growth parameter $\mu>0$ whose estimation is not straightforward. We solve this problem for a randomly generated dataset with $n=30000$ and $m=100$. We compare the following methods:
\begin{itemize}
\item FISTA \cite{Beck2009} with a fixed stepsize $\tau=\frac{1}{\hat L}$;
\item FISTA restart \cite{FISTArestart} with a fixed stepsize $\tau=\frac{1}{\hat L}$,
\item \texttt{FISTA\_adaBT} (\Cref{alg:FISTABT}) with $\rho=0.8$ and $\delta=0.95$,
\item Free-FISTA (\Cref{alg:FISTArestartBT_up}) with $\rho=0.8$ and $\delta=0.95$.
\end{itemize}

We set $\lambda_1=10$, $\lambda_2 = 3$ and $x_0\leftarrow\mathcal{U}\left([-1,1]\right)$. We get that $\hat L\geqslant 9\cdot10^5$ is an upper bound of $L$. An estimate of the solution of \eqref{eq:reglog} is pre-computed by running Free-FISTA for a large number of iterations. This allows us to compute for all methods $\log \left(F(r_j)-\hat F\right)$ with $\hat F\approx F^*$. 
\begin{figure}[!t]
\centering
\includegraphics[width=0.98\textwidth]{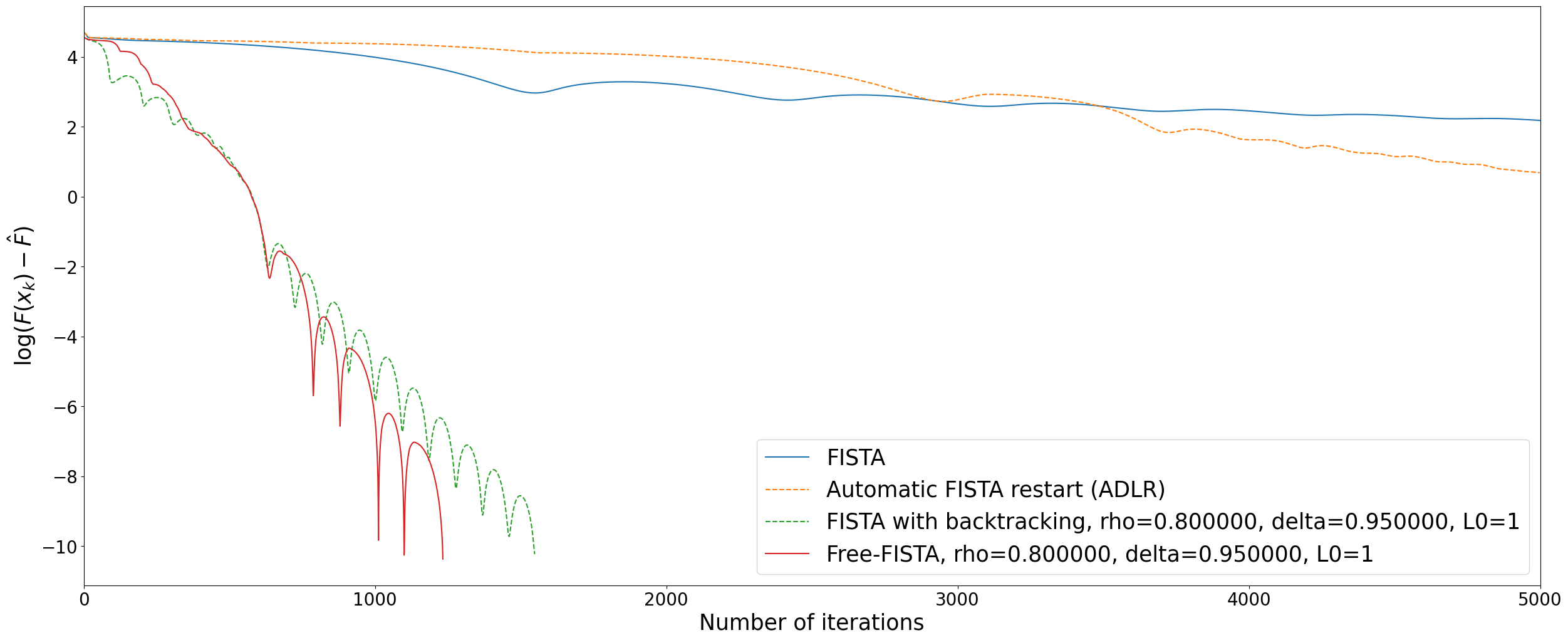}
\caption{Convergence rates w.r.t.~the number of total iterations (backtracking iterations are not taken in account) for problem \eqref{eq:reglog}.}
\label{fig:reglog_cvg}
\end{figure}

\begin{figure}[!t]
\centering
\includegraphics[width=0.98\textwidth]{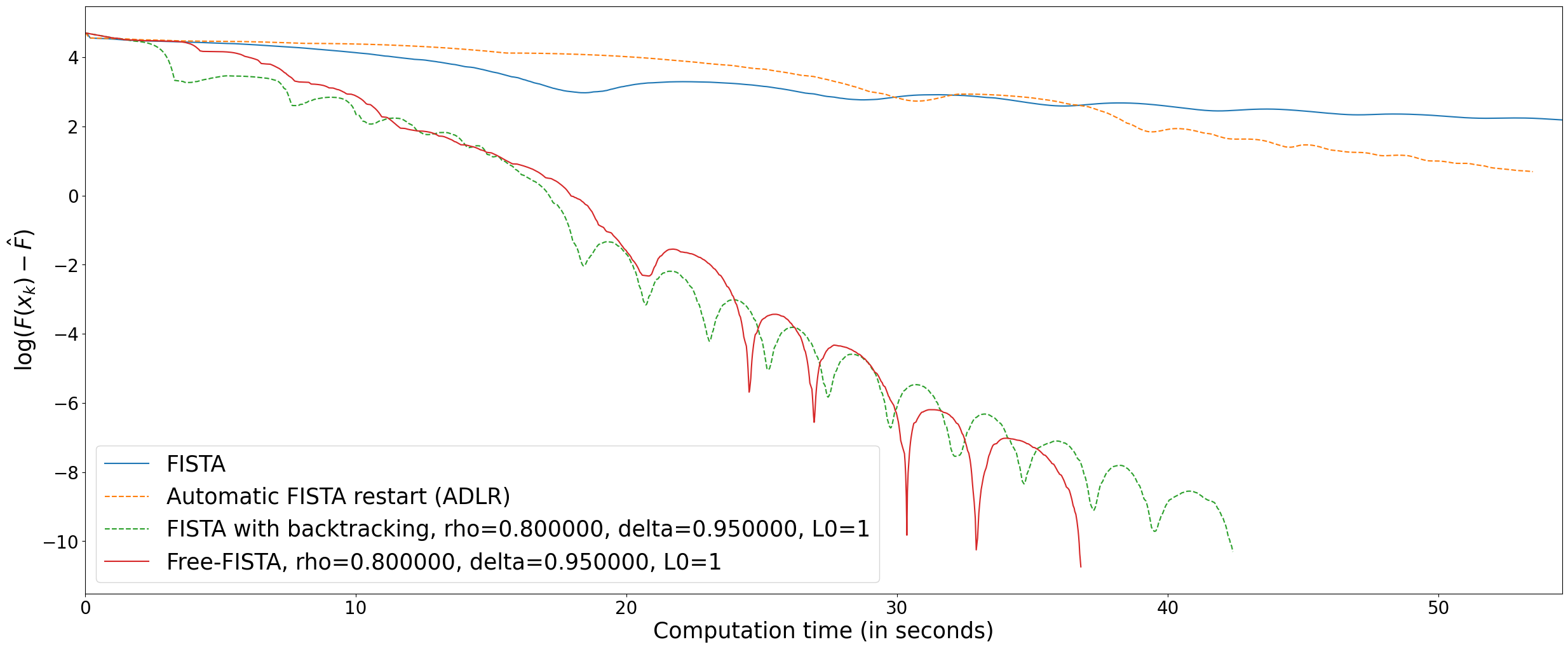}
\caption{Convergence rates w.r.t.~CPU times for problem \eqref{eq:reglog}.}
\label{fig:reglog_cvg_ctime}
\end{figure}

In Figure \ref{fig:reglog_cvg} the convergence rates of each algorithm are compared w.r.t.~the total number of iterations without taking into account the inner iterations required by the backtracking loops. We observe that the use of the adaptive backtracking accelerates both FISTA and FISTA restart. The improved efficiency provided by the combination of restarting and backtracking strategies is highlighted since Free-FISTA is the fastest method. Note, however, that an exhaustive information on the efficiency of each method can not directly be deduced by this plot as the computational burdens required by the use of the inner backtracking routines are not reported. We thus complement our considerations with Figure \ref{fig:reglog_cvg_ctime} which allows us to compare the methods w.r.t. the computation time. One can observe that the additional computations required by the backtracking strategy do not prevent the corresponding schemes from being faster.

\begin{figure}[!t]
\centering
\includegraphics[width=0.98\textwidth]{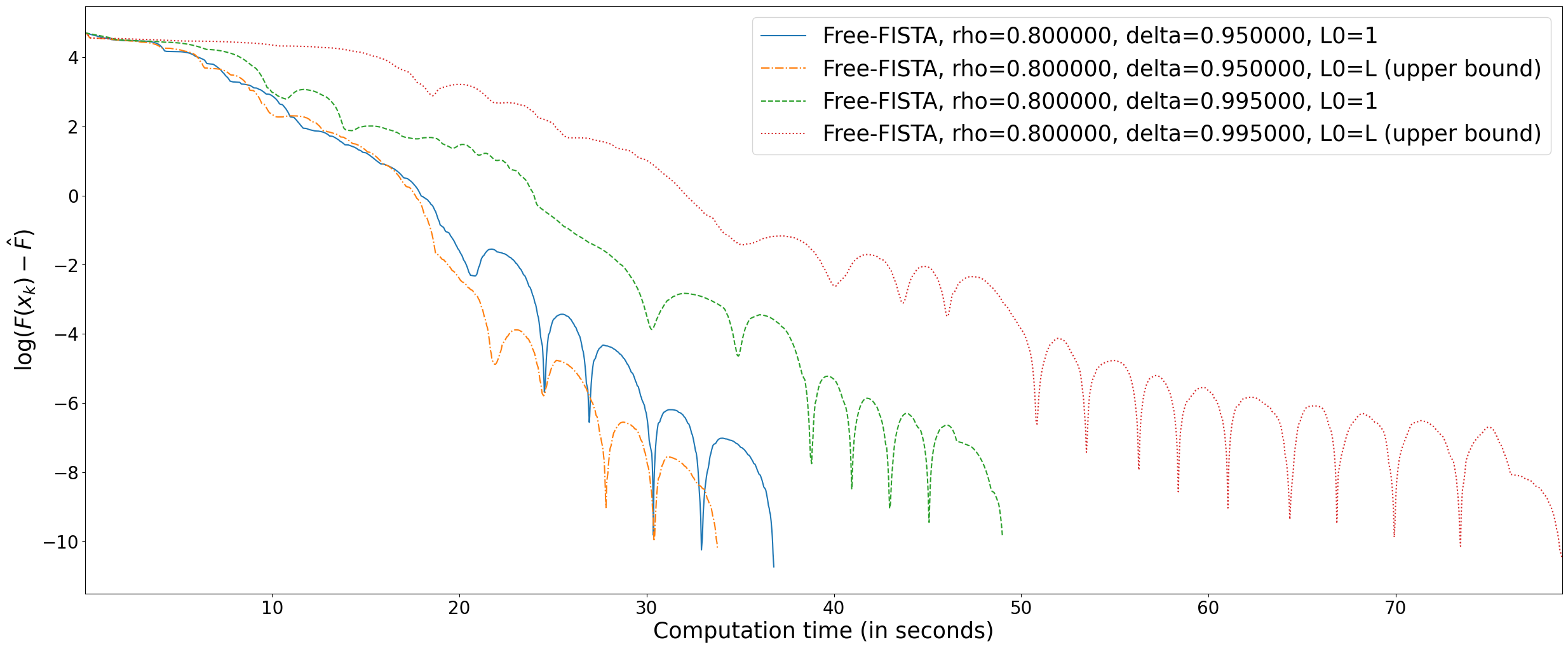}
\caption{Convergence rates of Free-FISTA for several choices of parameters $\rho$, $\delta$ and $L_0$ w.r.t. CPU time for problem \eqref{eq:reglog}.}
\label{fig:reglog_comp}
\end{figure}

Figure \ref{fig:reglog_comp} shows the convergence rate of Free-FISTA w.r.t. the computation time for several parameter choices. We take $\rho=0.8$, $\delta\in\{0.95,0.995\}$ and $L_0\in\{1,\hat L\}$ where $\hat L$ is the upper estimation of the Lipschitz constant of $\nabla f$ given in \eqref{eq:L_estimate_prob1} and $1$ is an arbitrary value. This graph shows that Free-FISTA is not highly sensitive to parameter variations in this example. Note that the choice $\delta=0.95$ seems to perform better than $\delta=0.995$. Indeed, as the Lipschitz constant of $\nabla f$ in this problem is poorly estimated, taking a small $\delta$ allows the scheme to explore different choices more efficiently. The value of $L_0$ has a small influence on the overall efficiency of the scheme.

\begin{figure}[!h]
\centering
\includegraphics[width=0.98\textwidth]{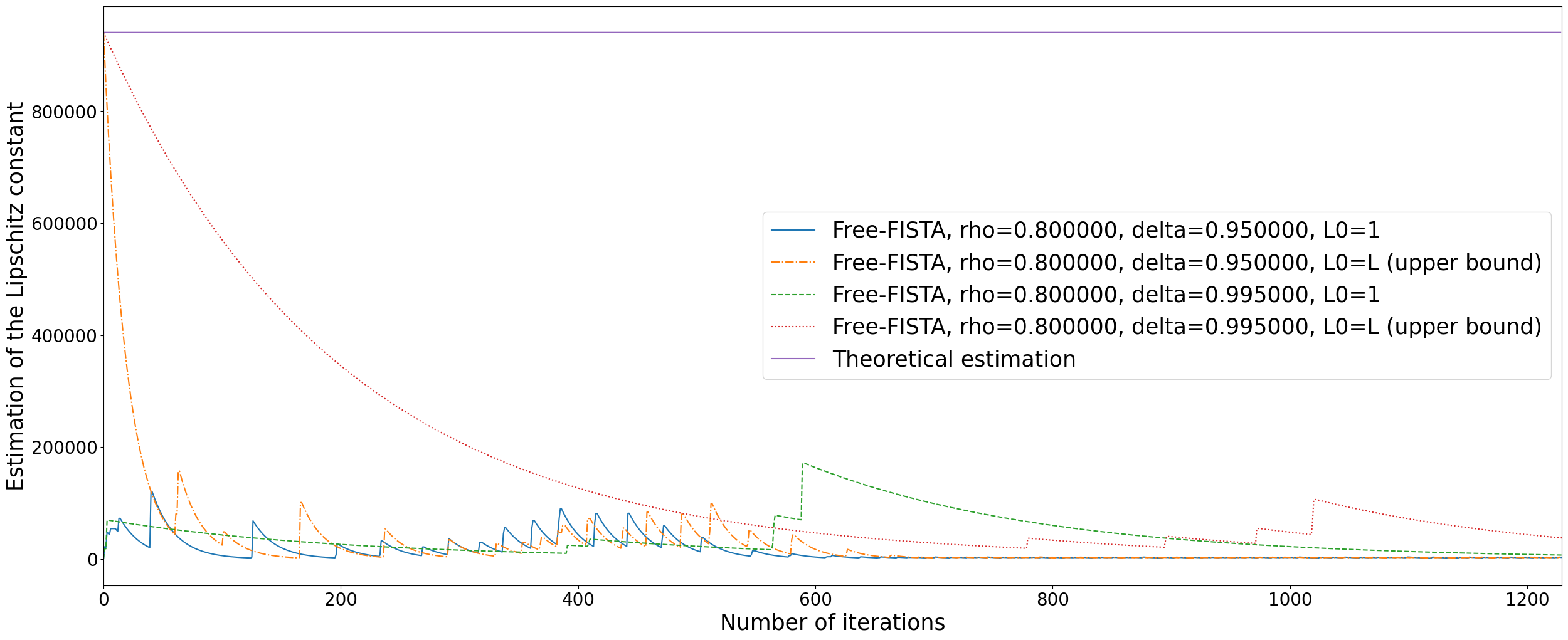}
\caption{Estimation of the Lipschitz constant of $\nabla f$ according to the number of FISTA iterations for problem \eqref{eq:reglog}.}
\label{fig:reglog_lip}
\end{figure}

Figure \ref{fig:reglog_lip} gives an overview of the estimations of the Lipschitz constant w.r.t. to FISTA iterations for each parameter choice. We can see that the theoretical upper bound $\hat L\geqslant 9\cdot 10^5$ is significantly large compared to the estimations computed by Free-FISTA for any set of parameters (the last estimates are approximately equal to $3000$). This explains the substantial performance gap between schemes involving a constant stepsize and backtracking methods (see Figure \ref{fig:reglog_cvg}) as a lower Lipschitz constant allows larger stepsizes. In addition, Figure \ref{fig:reglog_lip} shows that a lower value of $\delta$ encourages larger variations of estimates of $L$ per FISTA iteration, allowing for greater flexibilty. 

\begin{figure}[!h]
\centering
\includegraphics[width=0.98\textwidth]{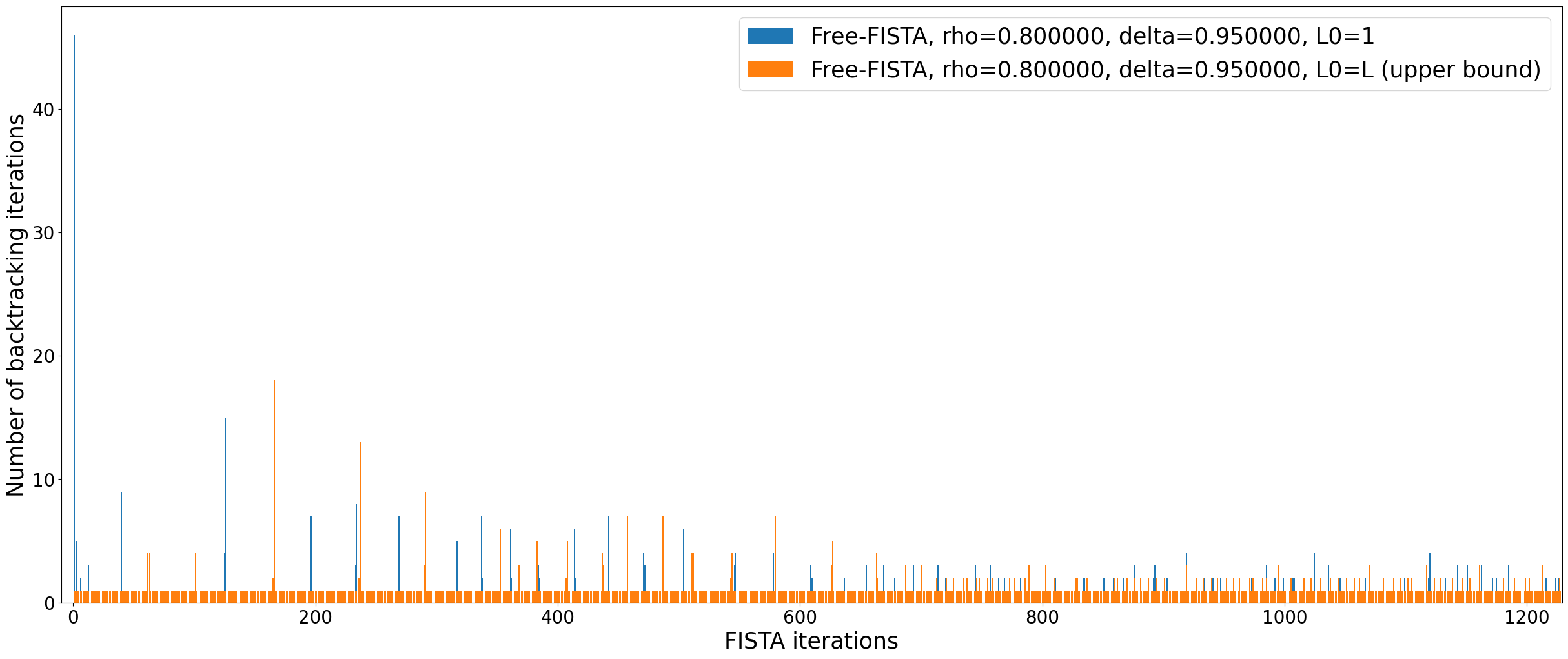}
\caption{Number of backtracking calls per total FISTA iterations for problem \eqref{eq:reglog}.}
\label{fig:reglog_ite}
\end{figure}

In Figure \ref{fig:reglog_ite}, we compare the differences observed between choosing a lower or an upper estimation $L_0$ of the Lipschitz constant $L$. Setting $L_0$ as a lower estimate forces the backtracking routine to compute a significant number of backtracking iterations before finding an estimate $\tilde L$ such that the stepsize $\frac{1}{\tilde L}$ is admissible. Once this is done, this estimation is generally tight and the number of backtracking iterations decreases critically. By taking $L_0$ as an upper estimate, we observe that the total number of backtracking iterations is smaller, but the estimation of $L$ stays poor for several Free-FISTA iterations (see Figure \ref{fig:reglog_lip}). Both approaches are equally efficient for this example because the high cost of the backtracking routines in the first case is compensated by the small stepsizes in the first FISTA iterations of the second case.

We now follow the experiments provided in \cite{fercoq2019adaptive} and consider the dataset \textit{dorothea} ($n=100000$ and $m=800$) with $\lambda_1=10$ and $\lambda_2=\frac{\lambda_1\|A^Tb\|^2}{80\|A^Tb\|_\infty n}=0.9097$. 
Table \ref{tab:reglog} compares the efficiency of the backtracking and restarting strategies for this example evaluated in terms of the CPU time required to satisfy the stopping condition with $\varepsilon=10^{-5}$. One can observe that methods involving adaptive backtracking are significantly faster. Algorithm \ref{alg:FISTArestartBT_up} is the most efficient algorithms, being, in addition, fully automatic. Some sensitivity to parameters $\rho$ and $\delta$ is observed, which, however, does not seem to significantly impact the overall computational gains.

\begin{table}[!h]
\centering
\small{
\begin{tabular}{|c | c | c | c|} 
 \hline
 \textbf{Algorithm} & $\rho$ & $\delta$ & Time (s) \\ [0.5ex] 
 \hline\hline
 FISTA & - & - & $28594$ \\ 
 \hline
 FISTA restart & - & - & $12825$ \\
  \hline
 \texttt{FISTA\_adaBT} & $0.85$ & $0.95$ & $3292$ \\
 & $0.8$ & $0.95$ & $2348$ \\
  \hline
 Free-FISTA & $0.85$ & $0.95$ & $1173$\\
 & $0.8$ & $0.95$ & $\mathbf{989}$\\[1ex]
 \hline
 \end{tabular}
 }
\caption{CPU times (mean over 5 runs) of different algorithms solving \eqref{eq:reglog} for the dataset \textit{dorothea} ($n=10^6$ and $m=800$), $\lambda_1=10$, $\lambda_2
=0.9097$ and $\varepsilon=10^{-5}$.}
\label{tab:reglog}
\end{table}

\subsection{Image inpainting}
\begin{figure}[!h]
\centering
\includegraphics[width=0.9\textwidth]{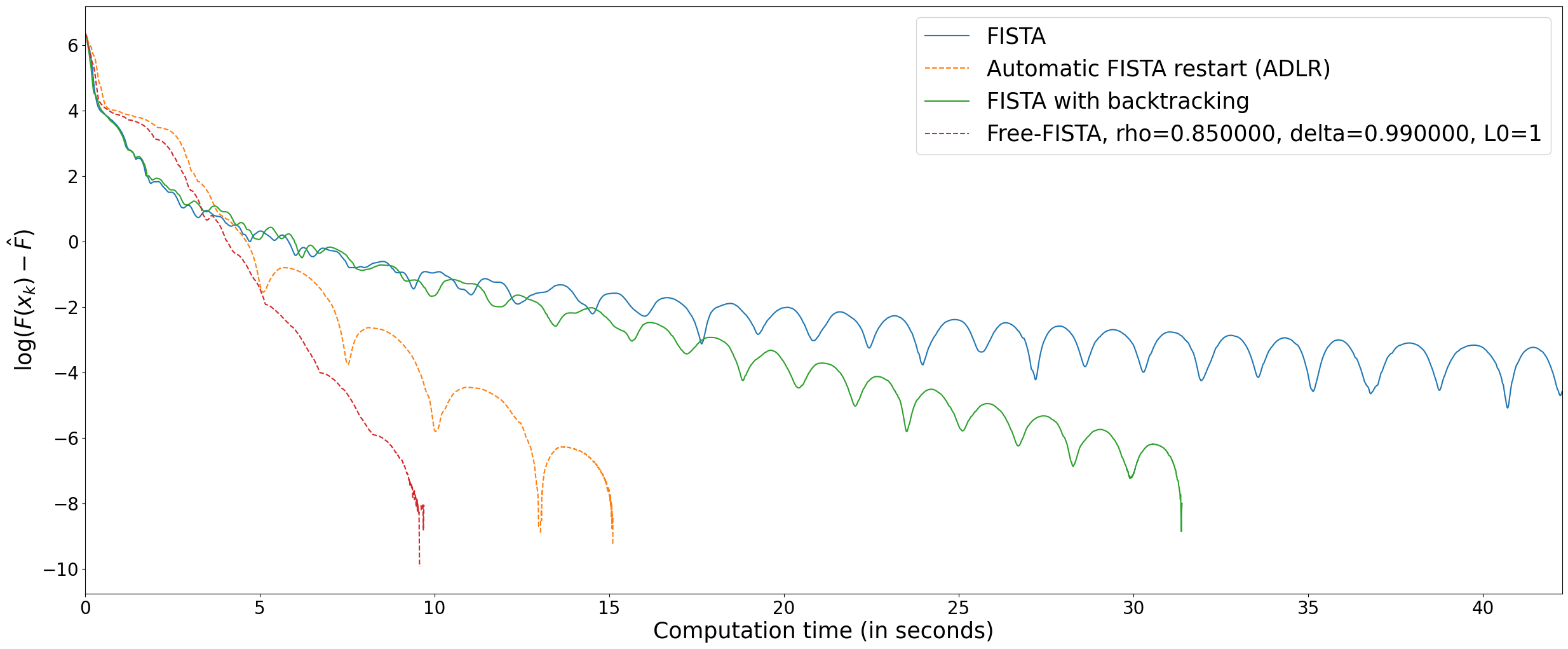}
\caption{Convergence rates in function values w.r.t the CPU time for problem \eqref{eq:inpainting}.}
\label{fig:inpainting_cvg_ctime}
\end{figure}

\begin{figure}[!h]
\centering
\includegraphics[width=0.9\textwidth]{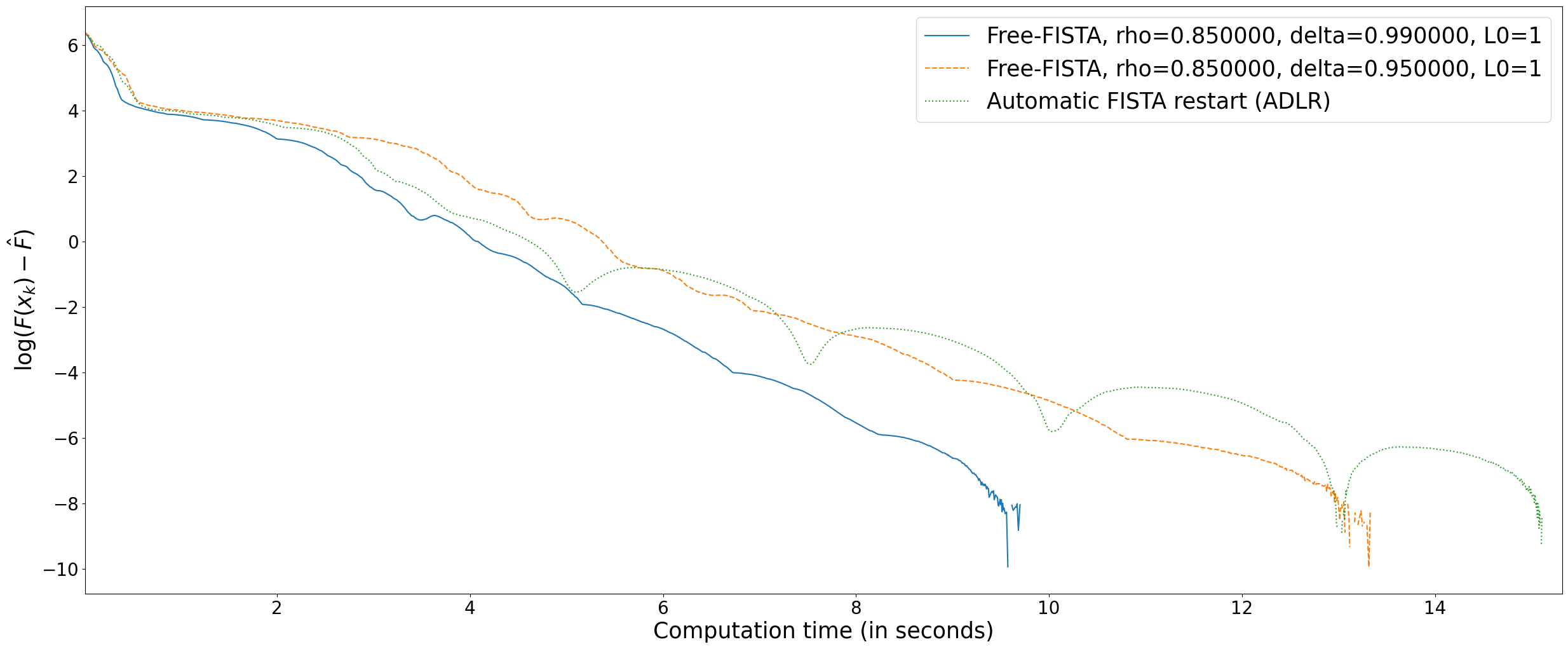}
\caption{Convergence rates in function valuesw.r.t the CPU time for problem \eqref{eq:inpainting}.}
\label{fig:inpainting_comp}
\end{figure}

\begin{figure}[!t]
\centering
\includegraphics[width=0.9\textwidth]{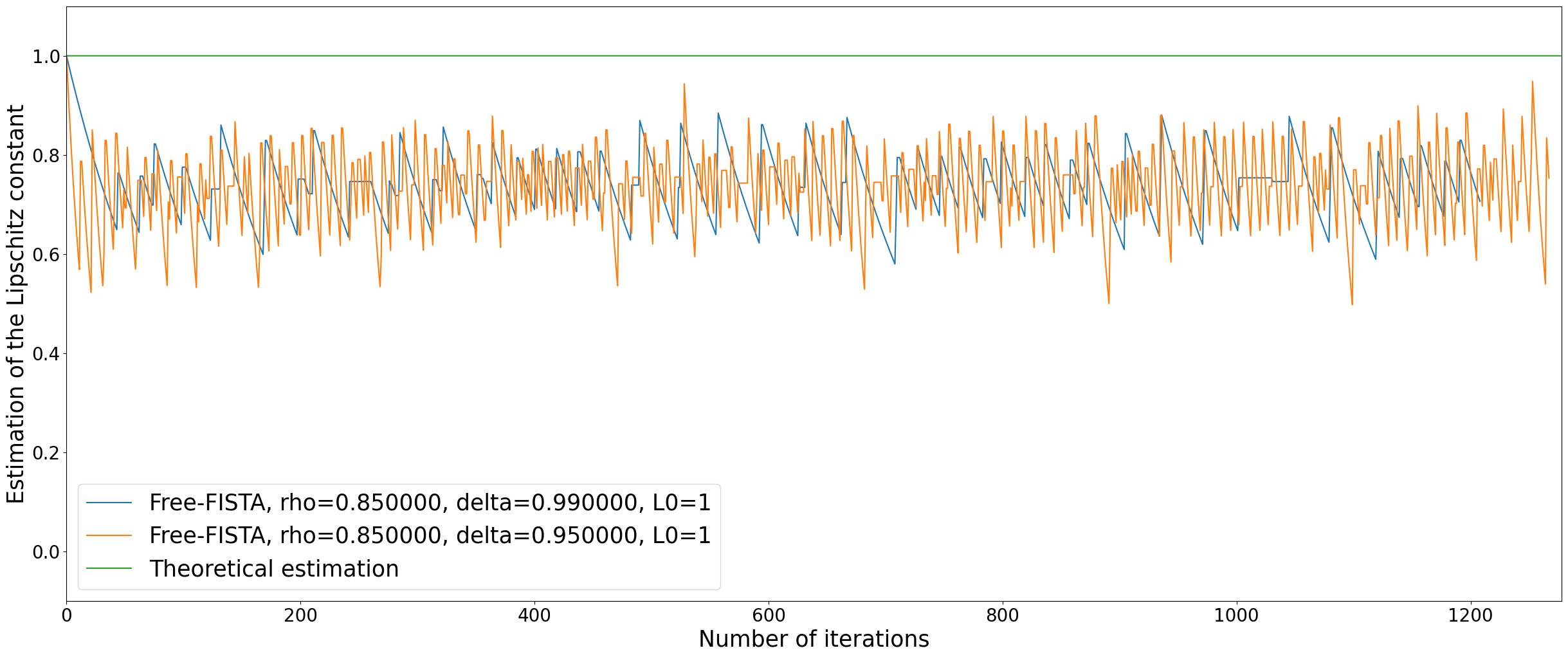}
\caption{Estimation of the Lipschitz constant of $\nabla f$ according to the number of FISTA iterations for problem \eqref{eq:inpainting}.}
\label{fig:inpainting_lip}
\end{figure}

\begin{figure}[h]
\centering
\includegraphics[width=0.9\textwidth]{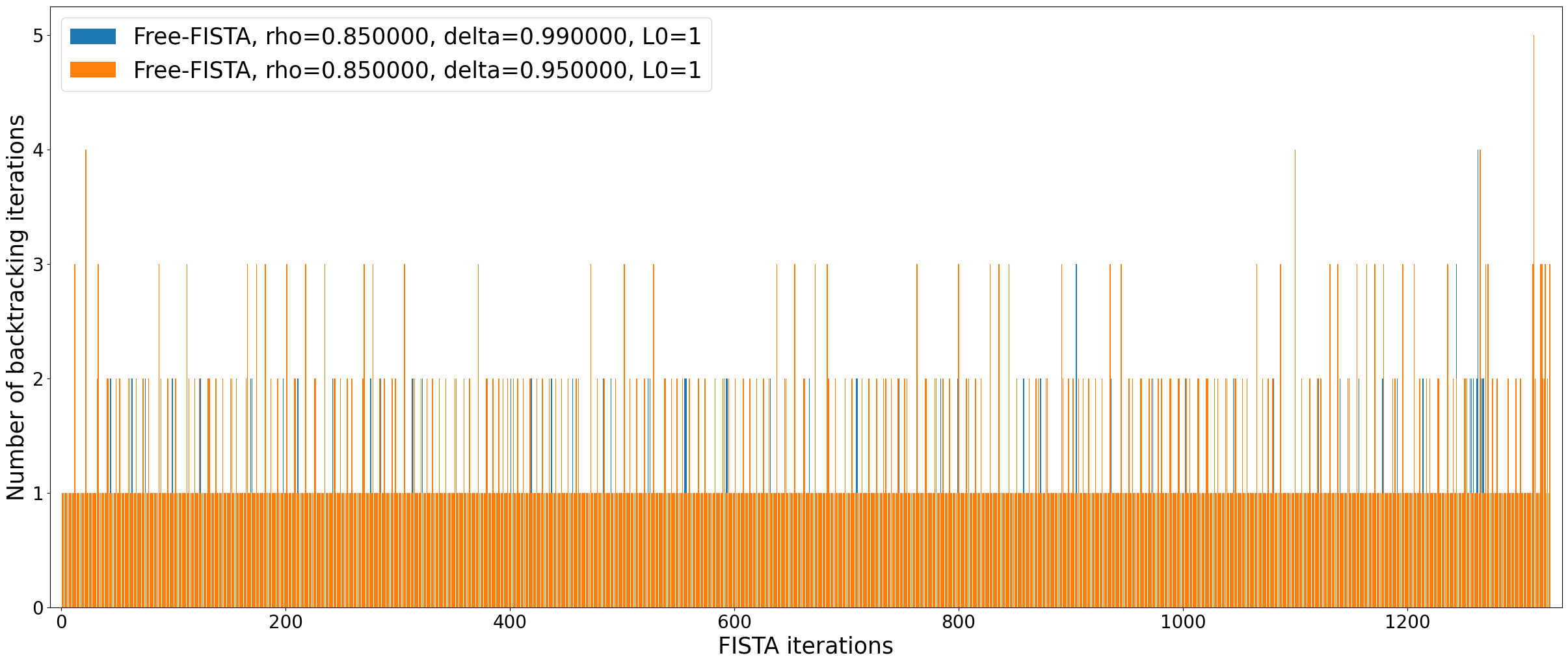}
\caption{Number of backtracking calls per total FISTA iterations for problem \eqref{eq:inpainting}.}
\label{fig:inpainting_ite}
\end{figure}

We now consider the problem of retrieving an image $\hat{x}\in\mathbb{R}^N$ from incomplete measurements $y=M\hat{x}$ where  $M\in\R^{N\times N}$ is a masking operator. We consider the regularized approach:
\begin{equation}
\argmin\limits_{x}~F(x):= f(x) + h(x) =\frac{1}{2}\|Mx-y\|^2+\lambda\|Tx\|_1,
\label{eq:inpainting}
\end{equation}
where $T\in\mathbb{R}^{N\times N}$ is an orthogonal transformation ensuring that $T\hat{x}$ is sparse. For this example we consider $\hat{x}$ to be piece-wise smooth, so that $T$ can be chosen  as an orthogonal wavelet transform. 
\begin{figure}[t!]
    \centering
\begin{minipage}{.48\textwidth}
\centering
   \includegraphics[height=3cm]{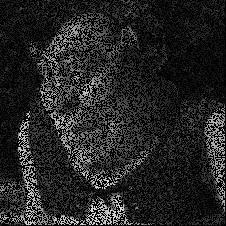}\quad\includegraphics[height=3cm]{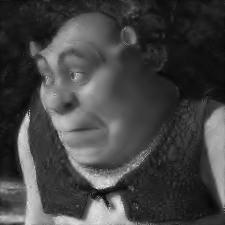}
\captionof{figure}{Data for problem \eqref{eq:inpainting}: the damaged image $y$ (left) and an inpainted result  (right).}
\label{fig:inpainting_ex}
\end{minipage}
\hfill 
\begin{minipage}{.48\textwidth}
\centering
   \includegraphics[height=3cm]{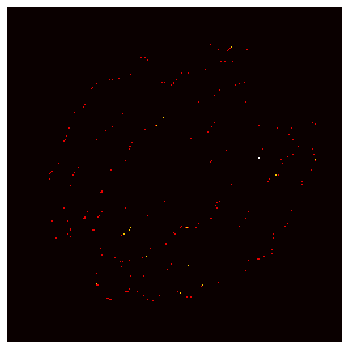}
\includegraphics[height=3
cm]{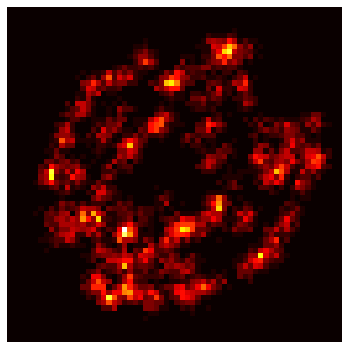}
\caption{Data for problem \eqref{eq:KL_l1}: ground-truth SMLM ISBI13 frame (left), and low-resolution data (right).}
\label{fig:deconv}
\end{minipage}
\end{figure}
%
%
The function $F$ satisfies the growth condition $\mathcal{G}^2_\mu$ for some $\mu>0$ which is not easily computable. In this case, it is trivial to show that an estimate of the Lipschitz constant of $\nabla f$ is $L=1$. Therefore, applying a backtracking strategy may seem superfluous as it involves additional computations. Nonetheless, we apply the methods previously introduced to test their performance with/without restarting. These tests are done on a picture with a resolution of $225\times225$ pixels, considering the wavelet Daubechies 4 and $\lambda=2$.
\begin{figure}[!h]
\centering
\includegraphics[width=0.9\textwidth]{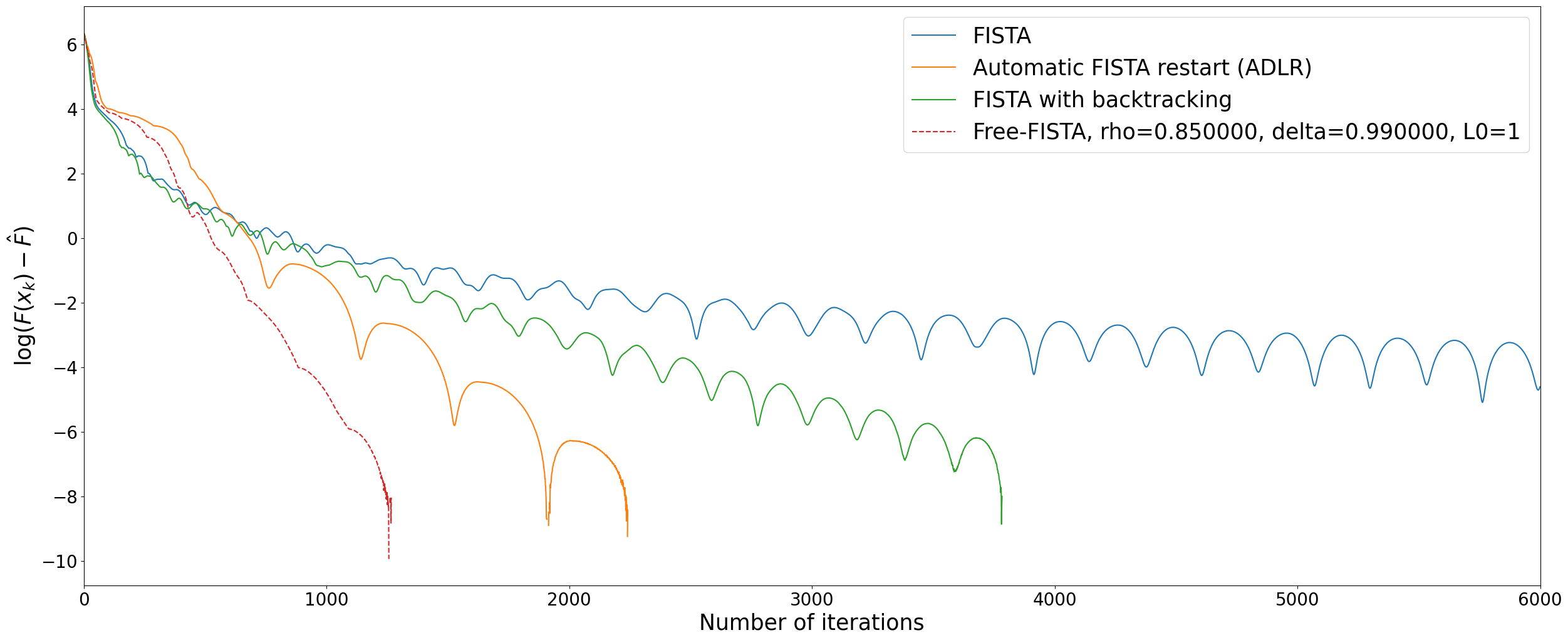}
\caption{Convergence rates in function values VS. total number of FISTA iterations (backtracking iterations are not taken in account) for problem \eqref{eq:inpainting}.}
\label{fig:inpainting_cvg}
\end{figure}
Figure \ref{fig:inpainting_cvg} shows that the backtracking procedure slightly improves the convergence of plain FISTA and FISTA restart w.r.t.~the total number of FISTA iterations. Observe that the benefits of backtracking are not as significant as in the previous example since the estimate of the Lipschitz constant $L=1$ is here accurate. In Figure \ref{fig:inpainting_cvg_ctime} we observe that the additional backtracking loops do not affect the efficiency of the schemes in terms of CPU time. In this example, evaluating $f$ is indeed not expensive which explains their low computational costs.
In Figure \ref{fig:inpainting_comp} we compare the performance of Free-FISTA for different values of $\delta$ and in comparison with ADLR. We observe that $\delta$ should be taken rather large in this case. Contrary to the previous example, if $\delta$ is small ($\delta=0.95$), Free-FISTA performs many unnecessary backtracking iterations to compensate for the over-estimation of the step-sizes, which results in longer CPU times. This can be observed in Figure \ref{fig:inpainting_lip} and Figure \ref{fig:inpainting_ite}. By taking $\delta=0.99$, a more gentle estimation with less variability of $L$ is observed over time, with fewer backtracking iterations per FISTA iteration.

\subsection{Poisson image super-resolution with $\ell^1$ regularization}

 As a last example, we consider the image super-resolution problem for images corrupted by Poisson noise, a problem encountered, for instance, in fluorescence microscopy applications\cite{Lazzaretti2021,RebegoldiCalatroni2022}.
Given a blurred and noisy image $z\in \R_+^m$, the problem consists in retrieving a sparse and non-negative image $x\in\R_+^n$ from $z=\mathcal{P}(MHx+b)\in\mathbb{R}^m$ with $m=q^2n$, $q>1$, where $M\in\mathbb{R}^{m\times n}$  is a $q$-down-sampling operator of factor, $H\in\R^{n\times n}$ is a convolution operator computed for a given point spread function (PSF), $b=\bar{b}e_m\in\R^m_{>0}$ is a positive constant background term\footnote{We use the notation $e_d$ to denote the vector of all ones in $\mathbb{R}^d$.} and $\mathcal{P}(w)$ denotes a realization of a Poisson-distributed $m$-dimensional random vector with parameter $w\in\R^m_+$. 
%
To model the presence of Poisson noise in the data, we consider the generalized Kullback-Leibler divergence functional \cite{Bertero2018} defined by:
\begin{equation}  \label{eq:KL}
    f(x) = KL(MHx+b;z):=\sum_{i=1}^m\left(z_i\log \frac{z_i}{(MHx)_i+ \bar{b}}+(MHx)_i+\bar{b}-z_i\right),
\end{equation}
and where  the convention $0\log 0=0$ is adopted. We enforce sparsity by means of a $\ell^1$ penalty and impose non-negativity of the solution using  the indicator function  $\iota_{\geq 0}(\cdot)$ of the non-negative orthant, so as to consider:
\begin{equation}  \label{eq:KL_l1}
\min_{x\in\R^n}~F(x):=KL(MHx+b,z)+\lambda\|x\|_1+i_{\geq 0}(x).
\end{equation}
%
We can compute $\nabla f(x) = (MH)^T e_m - (MH)^T \left( \frac{z}{MHx + b}\right).$
Following \cite{Harmany2012,RebegoldiCalatroni2022}, we have that
$\nabla f$ is Lipschitz continuous on $\left\{x: x\geq 0  \right\}$ and its Lipschitz constant $L$ can be overestimated by:
\begin{equation}  \label{eq:Poiss_L}
 L=\frac{\max z_i}{\bar{b}^2} \max ((MH)^T e_m) \max(MH e_n).
\end{equation}

\begin{figure}[H]
\centering
\includegraphics[width=0.96\textwidth]{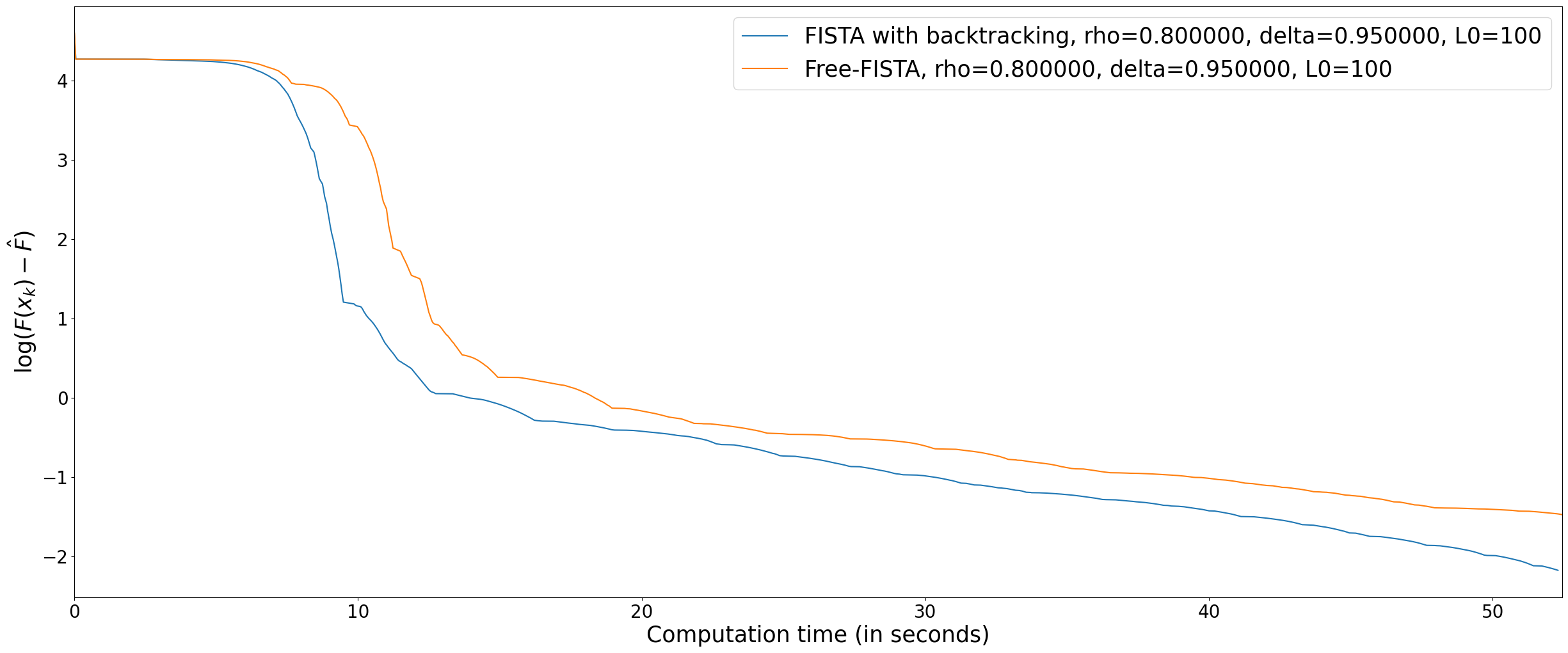}
\caption{Convergence rates in function values VS. CPU time for problem \eqref{eq:KL_l1}.}
\label{fig:deblurring_cvg}
\end{figure}

The theoretic estimation of $L$ in \eqref{eq:Poiss_L} may  be significantly large in particular, when $\bar{b}\ll 1$. Furthermore,
as showed in \cite{CalatroniI3D2021}, the Kullback-Leibler functional \eqref{eq:KL} is (locally) $2$-conditioned, hence $F$ satisfies $\mathcal{G}^2_\mu$ for some unknown $\mu>0$. The use of the Free-FISTA Algorithm \ref{alg:FISTArestartBT_up} thus seems appropriate. Results are showed in Figure \ref{fig:deblurring_cvg}.  For this problem, a clear advantage in the use of Free-FISTA in comparison with FISTA with adaptive backtracking cannot be observed. We observe that FISTA with adaptive backtracking is indeed faster in terms of iterations and consequently in terms of complexity (Free-FISTA requires additional computations being based on restarts). We argue that the inefficiency of the restarting strategy can be explained here by the geometry of $F$ in \eqref{eq:KL_l1}. The lack of any oscillatory behavior of FISTA endowed with adaptive backtracking suggests indeed that the function $F$ is flat, or, in other words, that $\mu$ is significantly small. Since restarting methods aim to handle the excess of inertia and oscillations, it appears not pertinent to apply such a method in this context. 

\appendix
\section{Proofs of the main results}
\label{sec:proof}

\subsection{Proof of \Cref{prop:Nesterov_up}}
\label{sec:prop1}
(i) As $F$ satisfies $\mathcal{H}_L$ for some $L>0$, \Cref{theo:convergence} combined with \eqref{eq:conv_FISTA_BT} states that the sequence $(x_k)_{k=1,\dots,n}$ provided by \Cref{alg:FISTABT} satisfies for all $k=1,\ldots,n$
\begin{equation*} 
F(x_{k+1})-F^*\leq \frac{2L}{\rho(k+1)^2} \|x_0-x^*\|^2,
\end{equation*}
 for all $x^*\in X^*$, whence
\begin{equation}
F(x_{k+1})-F^*\leqslant\frac{2L}{\rho(k+1)^2}~d(x_0,X^*)^2.\label{eq:FISTABT_bound2}
\end{equation}
Since $F$ further satisfies  $\mathcal{G}^2_\mu$ \eqref{eq:Lojasiewicz_dis}, we deduce \eqref{eq:Nesterov_up1} by combining \eqref{eq:Lojasiewicz_dis} and \eqref{eq:FISTABT_bound2}.

(ii) At each iteration $k\geq 0$ of \Cref{alg:FISTABT}, the following condition is satisfied:
\begin{equation*}
D_f(x_{k+1},y_{k+1})\leqslant\frac{\|x_{k+1}-y_{k+1}\|^2}{2\tau_{k+1}}.
\end{equation*}
As a consequence, the descent condition \eqref{condition:backtr} becomes:
\begin{align}
& F(x_{k+1})+\frac{\|x_{k+1}-x_k\|^2}{2\tau^{k+1}}\leqslant F(x_k)+\frac{\|y_{k+1}-x_k\|^2}{2\tau_{k+1}}  \notag\\
&\leqslant F(x_k)+\frac{(t_k-1)^2}{t_{k+1}^2}\frac{\|x_{k}-x_{k-1}\|^2}{2\tau_{k+1}} 
\leqslant F(x_k)+\frac{(t_k-1)^2}{t_{k+1}^2}\frac{\tau_k}{\tau_{k+1}}\frac{\|x_{k}-x_{k-1}\|^2}{2\tau_{k}}.  \label{eq:descent2}
\end{align}
By definition, there holds $\tau_{k+1}t_{k+1}(t_{k+1}-1)=\tau_kt_k^2$. Hence:
\begin{equation*}
\frac{(t_k-1)^2}{t_{k+1}^2}\frac{\tau_k}{\tau_{k+1}}=\frac{(t_k-1)^2t_{k+1}(t_{k+1}-1)}{t_k^2t_{k+1}^2}\leqslant1,
\end{equation*}
hence, from \eqref{eq:descent2} we get:
\[
F(x_{k+1})+\frac{\|x_{k+1}-x_k\|^2}{2\tau_{k+1}}\leqslant F(x_k)+\frac{\|x_{k}-x_{k-1}\|^2}{2\tau_{k}}\]
for all $k\geq 0$, whence we deduce \eqref{eq:Nesterov_up2}.


\subsection{Proof of \Cref{lem:mu_j}}
\label{sec:lem1}
Let $(\kappa_j)_{j\geqslant 2}$ be the sequence defined by
\begin{equation*}
\forall j\geqslant2,\quad \kappa_j:=\min\limits_{\underset{i<j}{i\in\N^*}}~\frac{4}{\rho(n_{i-1}+1)^2}\frac{F(r_{i-1})-F(r_{j})}{F(r_i)-F(r_{j})}.
\end{equation*}
We prove in this section that $(\kappa_j)_{j\geqslant2}$ is non increasing and bounded from below by the true inverse of the conditioning of the considered optimization problem.

First of all, according to \Cref{prop:Nesterov_up}, remember that we have \eqref{eq:ineq_muL} i.e.
$$
\forall i\in\N^*,\quad \kappa\leqslant\frac{4}{\rho(n_{i-1}+1)^2}\frac{F(r_{i-1})-F^*}{F(r_i)-F^*}.
$$
Since the application $u\mapsto \frac{F(r_{i-1})-u}{F(r_i)-u}$ is non decreasing on $[F^*,F(r_i)]$ (since $F(r_i)\leqslant F(r_{i-1})$), we deduce that for all $i\in \N^*$,
\begin{equation*}
\forall i<j,\quad \kappa\leqslant\frac{4}{\rho(n_{i-1}+1)^2}\frac{F(r_{i-1})-F(r_{j})}{F(r_i)-F(r_{j})}.
\end{equation*}
Hence, for a given $j\in \N^*$ and taking the infimum over the indexes $i\in \N^*$ such that $i<j$, we get: $\kappa \leqslant \kappa_j$. 
To complete the proof, we have that for all $j\geq 2$:
\begin{equation*}
\small{ \kappa_{j+1} =\min\limits_{\underset{i<j+1}{i\in\N^*}}~\frac{4}{\rho(n_{i-1}+1)^2}\frac{F(r_{i-1})-F(r_{j+1})}{F(r_i)-F(r_{j+1})}  \label{eq:monotone_kappa} \leqslant \min\limits_{\underset{i<j}{i\in\N^*}}~\frac{4}{\rho(n_{i-1}+1)^2}\frac{F(r_{i-1})-F(r_{j+1})}{F(r_i)-F(r_{j+1})} }
\end{equation*}
by simply observing that in \eqref{eq:monotone_kappa} the minimum is taken over a larger set.
By now applying \eqref{eq:Nesterov_up2} at the $j+1$ restart iteration we have that $F(r_{j+1})\leqslant F(r_j)$. As a consequence the function defined by $y\mapsto \frac{F(r_{i-1})-y}{F(r_i)-y}$ is an increasing homographic function which implies that for all $j\geq 2$:
\begin{equation*}
\small{
\kappa_{j+1}\leqslant \min\limits_{\underset{i<j}{i\in\N^*}}~\frac{4}{\rho(n_{i-1}+1)^2}\frac{F(r_{i-1})-F(r_{j+1})}{F(r_i)-F(r_{j+1})}\leqslant \min\limits_{\underset{i<j}{i\in\N^*}}~\frac{4}{\rho(n_{i-1}+1)^2}\frac{F(r_{i-1})-F(r_{j})}{F(r_i)-F(r_{j})}=  \kappa_j.}
\end{equation*}


\subsection{Proof of \Cref{lem:cor}}
\label{sec:lem4}

Suppose that $F$ satisfies $\mathcal{H}_L$ and $\mathcal{G}^2_\mu$ for some $L>0$ and $\mu>0$. Then, by \Cref{lem:Loja} 
\begin{equation*}
\forall x\in\R^N,\quad F(x)-F^*\leqslant \frac{2}{\mu}d(0,\partial F(x))^2.
\end{equation*}

Let now $x\in\R^N$ and $\tau>0$. By definition  \eqref{eq:prox_operator}, $x^+=T_{\tau} x$ is the unique minimizer of the function defined by $z\mapsto h(z)+\frac{1}{2\tau}\|z-x+\tau\nabla f(x)\|^2$. Thus, $T_\tau x$ satisfies
\begin{equation*}
0\in\partial h(T_\tau x)+\left\{\frac{1}{\tau}(T_\tau x-x)+\nabla f(x)\right\},
\end{equation*}
which entails: $g_\tau(x)-\nabla f(x)+\nabla f(T_\tau x)\in \partial F(T_\tau x)$. By the $L$-Lipschitz continuity of $\nabla f$ we can now deduce
\begin{align*}
& \|g_\tau(x)-\nabla f(x)+\nabla f(T_\tau x)\|\leqslant \|g_\tau(x)\|+\|\nabla f(T_\tau x)-\nabla f(x)\|\\
&\leqslant \|g_\tau(x)\|+L\|T_\tau x-x\| \leqslant (1+L\tau)\|g_\tau(x)\|.
\end{align*}
By combining all these inequalities we conclude that
\begin{equation*}
\small{
F(T_\tau x)-F^*\leqslant \frac{2}{\mu}d(0,\partial F(T_\tau x))^2\leqslant \frac{2}{\mu}\|g_\tau(x)-\nabla f(x)+\nabla f(T_\tau x)\|^2\leqslant \frac{2(1+L\tau)^2}{\mu}\|g_\tau(x)\|^2.}
\end{equation*}

\subsection{Sketch of the proof of \Cref{thm:1}}

Since the proof is rather technical, we split it into the following two parts:
\begin{enumerate}
  \item We show that there is at least one doubling step every $T$ iterations for a suitable $T$.  In particular:
  \begin{enumerate}
      \item We suppose that there is no doubling step from $j=s+1$ to $j=s+T$ for $s\geq 1$.
      \item We show a geometrical decrease of $(F(r_{j-1})-F(r_j))_{j\in\llbracket s+1,s+T\rrbracket}$ where the factor represents the gain of the $j$-th execution of \Cref{alg:FISTABT}.
      \item We state and apply \Cref{lem:grad} (whose proof is given in Subsection~\ref{sec:lem3}) to show that there exists an upper bound for $\|g_{1/{L}_{j-1}^+}(r_{j-1})\|$  depending on $F(r_{j-1})-F(r_j)$ for all $j\in\llbracket s+1,s+T\rrbracket$.
      \item We show that the geometrical decrease in (b) entails that the exit condition $\|g_{1/{L}_{j-1}^+}(r_{j-1})\|\leq \varepsilon$ is satisfied for $j=s+T$.
  \end{enumerate}
  \item We use 1. to show that the total number of restarting iterations $\sum_{i=0}^jn_i$ is necessarily bounded by $2Tn_j$. The conclusion of \Cref{thm:1} thus comes from \Cref{lem:n_j} providing an upper bound of $n_j$.
\end{enumerate}

\subsection{Proof of \Cref{thm:1}}  \label{sec:proof_thm1}
Let $C>\frac{4}{\sqrt{\rho}}$ and $\varepsilon>0$. We first define 
\[
T:=1+\left\lceil\frac{\log\left(1+\frac{16}{C^2\rho-16}\frac{2L(F(r_0)-F^*)}{\rho\varepsilon^2}\right)}{\log\left(\frac{C^2\rho}{4}-1\right)}\right\rceil.
\]
We claim that a doubling step is performed at least every $T$ iterations.\\
For $s\geqslant1$, assume that there is no doubling step for $T-1$ iterations from $j=s+1$ to $j=s+T$. This means:
\begin{equation}\forall j\in\llbracket s+1,s+T\rrbracket ,\quad n_{j-1}>C\sqrt{\frac{1}{ \kappa_j}},\label{eq:no_doubling1_up}
\end{equation}
whence:
\begin{equation}\forall j\in\llbracket s,s+T\rrbracket ,\quad n_j=n_s,\label{eq:no_doubling2_up}\end{equation}
where the case $j=s$ trivially holds.
We deduce that $\forall j\in\llbracket s+2,s+T\rrbracket$:
\begin{align*}
{\kappa}_j&=\min\limits_{\underset{i<j}{i\in\N^*}}\frac{4}{\rho(n_{i-1}+1)^2}\frac{F(r_{i-1})-F(r_{j})}{F(r_i)-F(r_{j})}\leqslant\min\limits_{\underset{s<i<j}{i\in\N^*}}\frac{4}{\rho(n_{i-1}+1)^2}\frac{F(r_{i-1})-F(r_{j})}{F(r_i)-F(r_{j})}\\
&\leqslant\min\limits_{\underset{s<i<j}{i\in\N^*}}\frac{4}{\rho{n_{i-1}}^2}\frac{F(r_{i-1})-F(r_{j})}{F(r_i)-F(r_{j})}\leqslant\min\limits_{\underset{s<i<j}{i\in\N^*}}\frac{4}{\rho{n_{s}}^2}\frac{F(r_{i-1})-F(r_{j})}{F(r_i)-F(r_{j})}\\
&\leqslant\frac{4}{\rho{n_{s}}^2}\min\limits_{\underset{s<i<j}{i\in\N^*}}\frac{F(r_{i-1})-F(r_{j})}{F(r_i)-F(r_{j})},
\end{align*}
due to \eqref{eq:no_doubling2_up}. Using \eqref{eq:Nesterov_up2}, we deduce that:
\begin{equation}\forall j\in\llbracket s+2,s+T\rrbracket,\quad {\kappa}_j\leqslant\frac{4}{\rho{n_{s}}^2}\frac{F(r_{j-2})-F(r_{j})}{F(r_{j-1})-F(r_{j})}.\label{eq:mu_j1_up}\end{equation}
Combining now \eqref{eq:no_doubling1_up} with \eqref{eq:no_doubling2_up} and \eqref{eq:mu_j1_up} we get:
\begin{equation*}
n_s>C\sqrt{\frac{1}{\frac{4}{\rho{n_{s}}^2}\frac{F(r_{j-2})-F(r_{j})}{F(r_{j-1})-F(r_{j})}}}
=n_s\frac{C\sqrt{\rho}}{2}\sqrt{\frac{F(r_{j-1})-F(r_{j})}{F(r_{j-2})-F(r_{j})}}
\end{equation*}
which leads to
\begin{equation*}F(r_{j-2})-F(r_{j})>\frac{C^2\rho}{4}(F(r_{j-1})-F(r_{j})),\end{equation*}
which further entails
\begin{equation*}F(r_{j-2})-F(r_{j-1})>\left(\frac{C^2\rho}{4}-1\right)(F(r_{j-1})-F(r_{j})).\end{equation*}
Since $C>\frac{4}{\sqrt{\rho}}>\frac{2}{\sqrt{\rho}}$ we now get the following geometric functional decrease.
\begin{equation}F(r_{j-1})-F(r_j)<\frac{4}{C^2\rho-4}(F(r_{j-2})-F(r_{j-1})).\label{eq:ineq2_up}\end{equation}
We now consider the case $j=s+1$:
\begin{align*}
 \kappa_{s+1}&=\min\limits_{\underset{i<s+1}{i\in\N^*}}\frac{4}{\rho(n_{i-1}+1)^2}\frac{F(r_{i-1})-F(r_{s+1})}{F(r_i)-F(r_{s+1})}\leqslant \frac{4}{\rho(n_{s-1}+1)^2}\frac{F(r_{s-1})-F(r_{s+1})}{F(r_s)-F(r_{s+1})}\\
&\leqslant \frac{4}{\rho(\frac{n_s}{2}+1)^2}\frac{F(r_{s-1})-F(r_{s+1})}{F(r_s)-F(r_{s+1})}
\leqslant \frac{16}{\rho{n_s}^2}\frac{F(r_{s-1})-F(r_{s+1})}{F(r_s)-F(r_{s+1})},\notag
\end{align*}
since $n_s\leqslant 2n_{s-1}$.
By reapplying $C>\frac{4}{\sqrt{\rho}}$, similar computations show that
\begin{equation}
F(r_{s})-F(r_{s+1})<\frac{16}{C^2\rho-16}(F(r_{s-1})-F(r_{s})).
\label{eq:ineq1_up}
\end{equation}
To carry on with the proof, we now state \Cref{lem:grad} which links the composite gradient mapping $g$ to the function $F$. The proof is reported in \Cref{sec:lem3}:
\begin{lemma} \label{lem:grad}
Let $F$ satisfy the assumption $\mathcal{H}_L$ for some $L>0$. Then the sequence $(r_j)_{j\in\N}$ provided by \Cref{alg:FISTArestartBT_up} satisfies
\begin{equation*}
\forall j \geqslant1,\quad \frac{\rho}{2L}\|g_{1/{L}_{j}^+}(r_{j})\|^2\leqslant F(r_{j})-F(r_{j+1}),
\end{equation*}
where ${L}_{j}^+$ is an estimate of $L$ provided by \Cref{alg:FB_BT_up}.
\end{lemma}

By  \Cref{lem:grad} and recalling inequalities \eqref{eq:ineq2_up} and \eqref{eq:ineq1_up}, we can thus obtain the following sequence of inequalities
\begin{align*}
&\frac{\rho}{2L}\|g_{1/{L}_{s+T-1}^+}(r_{s+T-1})\|^2\leqslant F(r_{s+T-1})-F(r_{s+T}) \\
&\leqslant\frac{4}{C^2\rho-4}(F(r_{s+T-2})-F(r_{s+T-1}))\\
&\leqslant\left(\frac{4}{C^2\rho-4}\right)^{T-1}\left(\frac{16}{C^2\rho-16}\right)(F(r_{s-1})-F(r_{s}))\\
& \leqslant\left(\frac{4}{C^2\rho-4}\right)^{T-1}\left(\frac{16}{C^2\rho-16}\right)(F(r_0)-F^*)\\
&\leqslant\left(\frac{4}{C^2\rho-4}\right)^{\left\lceil\frac{\log\left(1+\frac{16}{C^2\rho-16}\frac{2L(F(r_0)-F^*)}{\rho\varepsilon^2}\right)}{\log\left(\frac{C^2\rho}{4}-1\right)}\right\rceil}\left(\frac{16}{C^2\rho-16}\right)(F(r_0)-F^*)\\
&\leqslant\left(\frac{4}{C^2\rho-4}\right)^{\frac{\log\left(1+\frac{16}{C^2\rho-16}\frac{2L(F(r_0)-F^*)}{\rho\varepsilon^2}\right)}{\log\left(\frac{C^2\rho}{4}-1\right)}}\left(\frac{16}{C^2\rho-16}\right)(F(r_0)-F^*)\\
&\leqslant\frac{1}{1+\frac{16}{C^2\rho-16}\frac{2L(F(r_0)-F^*)}{\rho\varepsilon^2}}\left(\frac{16}{C^2\rho-16}\right)(F(r_0)-F^*)\leqslant\frac{\rho\varepsilon^2}{2L}.
\end{align*}
As a consequence, if there are $T$ consecutive restarts without any doubling of the number of iterations, then the exit condition $\|g_{1/{L}_j^+}(r_j)\|\leqslant\varepsilon$ is eventually satisfied. This means that there exists a doubling step at least every $T$ steps and that for all $s\geqslant1$ there exists $j\in\llbracket s+1,s+T\rrbracket$ such that
\begin{equation*}
n_{j-1}< C\sqrt{\frac{1}{\kappa_j}},
\end{equation*}
which implies that $n_j=2n_{j-1}$. Now, since $(n_j)_{j\in\N}$ is an increasing sequence, we get that $n_{s+T}\geqslant n_j=2n_{j-1}\geqslant2n_s$, so that
\begin{equation}
n_s\leqslant\frac{n_{s+T}}{2},~\quad \forall s\geqslant 1.
\label{eq:double}
\end{equation}
Let us now rewrite $j$ as $j=m+nT$ where $0\leqslant m<T$ and $n\geqslant0$. By monotonicity of $(n_j)_{j\in\N}$ we have
\begin{equation*}
\sum_{i=0}^jn_i=\sum_{i=0}^{m+nT}n_i=\sum_{i=0}^mn_i+\sum_{l=0}^{n-1}\sum_{i=1}^Tn_{m+i+lT} \leqslant
T\sum_{l=0}^n n_{m+lT}=T\sum_{l=0}^n n_{j-lT}. 
\end{equation*}
According to equation \eqref{eq:double} we have $n_{j-T}\leqslant\frac{n_j}{2}$, that is
\begin{equation*}
n_{j-lT}\leqslant\left(\frac{1}{2}\right)^ln_j,\quad \forall l\in[\![0,n]\!].
\end{equation*}
We thus obtain the following inequalities
\begin{equation}\sum_{i=0}^jn_i\leqslant T\sum_{l=0}^n n_{j-lT}\leqslant T\sum_{l=0}^n \left(\frac{1}{2}\right)^ln_j\leqslant T\sum_{l=0}^\infty \left(\frac{1}{2}\right)^ln_j=2Tn_j.\label{eq:sommes}\end{equation}
Combining \eqref{eq:sommes} with \Cref{lem:n_j} we thus finally get the desired result for $j>0$
\small
\begin{align*}
\sum_{i=0}^j n_i&\leqslant2Tn_j\leqslant4C\sqrt{\frac{L}{\mu}}T\leqslant4C\sqrt{\frac{L}{\mu}}\left(1+\left\lceil\frac{\log\left(1+\frac{16}{C^2\rho-16}\frac{2L(F(r_0)-F^*)}{\rho\varepsilon^2}\right)}{\log\left(\frac{C^2\rho}{4}-1\right)}\right\rceil\right)\\
&\leqslant \frac{4C}{\log\left(\frac{C^2\rho}{4}-1\right)}\sqrt{\frac{L}{\mu}}\left(2\log\left(\frac{C^2\rho}{4}-1\right)+\log\left(1+\frac{16}{C^2\rho-16}\frac{2L(F(r_0)-F^*)}{\rho\varepsilon^2}\right)\right).
\end{align*}\normalsize

\subsection{Proof of \Cref{cor:1}}  \label{sec:proofcor1}
Let $F$ satisfy $\mathcal{H}_L$ and $\mathcal{G}^2_\mu$ for some $L>0$ and $\mu>0$. Let $(r_j)_{j\in\N}$ and $(n_j)_{j\in\N}$ be the sequences provided by \Cref{alg:FISTArestartBT_up} with $C>4/\sqrt{\rho}$, $\varepsilon>0$ and let $L_{min}\in(0,L)$. 
We consider the case where the exit condition $\|g_{1/{L}_j^+}(r_j)\|\leqslant \varepsilon$ is satisfied at first for at least $8C\sqrt{\frac{1}{\kappa}}$ iterations. We define the function $\psi_\mu:\R_+^*\rightarrow \left(8C\sqrt{\frac{1}{\kappa}},+\infty\right)$ by:
\small
\begin{equation*}
\psi_\mu:\gamma\mapsto\frac{4C}{\log\left(\frac{C^2\rho}{4}-1\right)}\sqrt{\frac{L}{\mu}}\left(2\log\left(\frac{C^2\rho}{4}-1\right)+\log\left(1+\frac{16}{C^2\rho-16}\frac{2L(F(r_0)-F^*)}{\rho\gamma}\right)\right).
\end{equation*}
\normalsize
By \Cref{thm:1}, the number of iterations required to ensure $\|g_{1/{L}_j^+}(r_j)\|\leqslant \varepsilon$ satisfies
$\sum_{i=0}^jn_i\leqslant \psi_\mu(\varepsilon^2)$.
As $\psi_\mu$ is strictly decreasing and $\sum_{i=0}^jn_i>8C\sqrt{\frac{L}{\mu}}$, we deduce:
\begin{equation*}
\psi_\mu^{-1}\left(\sum_{i=0}^jn_i\right)\geqslant\varepsilon^2,
\end{equation*}
where $\psi_\mu^{-1}$ is the inverse function of $\psi_\mu$. By now applying \Cref{lem:cor} and since by construction ${L}_j^+\geqslant L_{min}$, we get:
\begin{align}
\small{
F({r}_j^+)-F^*\leqslant\frac{2\left(1+\frac{L}{{L}_j^+}\right)^2}{\mu}\|g_{1/{L}_j^+}(r_j)\|^2
\leqslant \frac{2\left(1+\frac{L}{L_{min}}\right)^2}{\mu}\psi_\mu^{-1}\left(\sum_{i=0}^jn_i\right).  \label{theo:ineq11}}
\end{align}
Elementary computations show that:
\begin{equation*}
\psi_\mu^{-1}:n\mapsto \frac{2L}{\rho}\frac{16}{C^2\rho-16}\frac{1}{e^{-2\log(\frac{C^2\rho}{4}-1)}e^{\frac{\log(\frac{C^2\rho}{4}-1)}{4C}\sqrt{\frac{\mu}{L}}n}-1}(F(r_0)-F^*),
\end{equation*}
hence from \eqref{theo:ineq11}, we get:
\small
\begin{equation*}
F({r}_j^+)-F^* \leqslant \frac{4L\left(1+\frac{L}{L_{min}}\right)^2}{\rho\mu}\frac{16}{C^2\rho-16}\frac{1}{e^{-2\log(\frac{C^2\rho}{4}-1)}e^{\frac{\log(\frac{C^2\rho}{4}-1)}{4C}\sqrt{\frac{\mu}{L}}\sum_{i=0}^jn_i}-1}(F(r_0)-F^*)
\end{equation*}
\normalsize
We can thus conclude that
\begin{equation*}
F({r}_j^+)-F^*= \mathcal{O}\left(e^{-\frac{\log(\frac{C^2\rho}{4}-1)}{4C}\sqrt{\kappa}\sum_{i=0}^jn_i}\right).
\end{equation*}
We can further maximize the function $C\mapsto\frac{\log(\frac{C^2\rho}{4}-1)}{4C}$ to obtain the optimal value $\hat C\approx6.38/\sqrt{\rho}$. This choice leads to the desired convergence rate:
\begin{equation*}
F({r}_j^+)-F^*= \mathcal{O}\left(e^{-\frac{\sqrt{\rho}}{12}\sqrt{\kappa}\sum_{i=0}^jn_i}\right).
\end{equation*}

To conclude the proof, let now $\left(x_{k,j}\right)_{k\in\llbracket0,n_j\rrbracket}$ and $\left(\tau_{k,j}\right)_{k\in\llbracket0,n_j\rrbracket}$ denote the iterates of \Cref{alg:FISTABT} following the $j$-th restart and the corresponding step-sizes, respectively. Note that in particular we have $x_{0,j}= r_{j-1}^+$ and $x_{n_j,j}=r_j$. By applying standard arguments as in the proof of \Cref{prop:Nesterov_up} (see Section \ref{sec:prop1}) we deduce that for any $j\geqslant0$ and every $k>0$:
\begin{equation*}
    F(x_{k,j})+\frac{\|x_{k,j}-x_{k-1,j}\|^2}{2\tau_{k,j}}\leqslant F(x_{0,j}).
\end{equation*}
Such inequality thus entails:
\begin{equation*}
\|x_{k,j}-x_{k-1,j}\|^2\leqslant 2\tau_{k,j}\left(F( r_j^+)-F^*\right)\leqslant \frac{2}{L_{min}}\left(F(r_j^+)-F^*\right).
\end{equation*}
By applying the first claim of this Corollary on the right hand side of the inequality above, we guarantee the existence of $M>0$ such that for $j$ large enough:
\begin{equation*}
    \forall k\in\llbracket1,n_j\rrbracket,\quad \|x_{k,j}-x_{k-1,j}\|^2\leqslant \frac{2M}{L_{min}}e^{-\frac{\log(\frac{C^2\rho}{4}-1)}{4C}\sqrt{\kappa}\sum_{i=0}^jn_i},
\end{equation*}
which implies that $
\sum_{j,k}\|x_{k,j}-x_{k-1,j}\|<+\infty,
$
showing that the trajectory of the total number of FISTA iterates has finite length.

\subsection{Proof of \Cref{lem:grad}}
\label{sec:lem3}
Since by definition $({r}_j^+,{L}_j^+)=\text{\texttt{FB\_BT}}(r_j,{L}_j;\rho)$, for all $j\geqslant1$ there holds:
$
D_f({r}_j^+,r_j)\leqslant \frac{{L}_j^+}{2}\|{r}_j^+-r_j\|^2,$
with ${r}_j^+=T_{1/{L}_j^+}(r_j)$ which allows us to specialize the descent condition \eqref{condition:backtr} as:
\begin{equation*}
F({r}_j^+) + \frac{{L}_j^+}{2}\|{r}_j^+-x\|^2\leqslant F(x)+\frac{{L}_j^+}{2}\|r_j-x\|^2,
\end{equation*}
for all $x\in\R^N$. By choosing $x=r_j$ and by definition of $g_{1/{L}_j^+}$ we get:
\begin{equation*}
\frac{1}{2L_j^+}\|g_{1/L_j^+}r_j\|^2\leqslant F(r_j)-F({r}_j^+).
\end{equation*}
Since by \eqref{eq:Lk_rho}, we further deduce $L_j^+\leqslant\frac{L}{\rho}$,
\begin{equation*}
\frac{\rho}{2L}\|g_{1/L_j^+}(r_j)\|^2\leqslant F(r_j)-F({r}_j^+).
\end{equation*}
Inequality \eqref{eq:Nesterov_up2} ensures $F(r_{j+1})\leqslant F({r}_j^+)$ which finally entails.
\begin{equation*}
\frac{\rho}{2L}\|g_{1/L_j^+}(r_j)\|^2\leqslant F(r_j)-F(r_{j+1}).
\end{equation*}

\section*{Acknowledgements}

JFA and LC acknowledge the  support of the  EU  Horizon  2020 research and innovation program under the Marie Sk\l odowska-Curie NoMADS grant agreement No777826. LC acknowledges the support  of the ANR TASKABILE (ANR-22-CE48-0010). This work was supported by the ANR MICROBLIND (grant ANR-21-CE48-0008), the ANR Masdol (grant ANR-PRC-CE23), the FMJH Program PGMO 2019-0024 and EDF-Thales-Orange.

\bibliographystyle{siamplain}
\bibliography{FISTAbib}
\end{document}